\documentclass{amsart}
\usepackage{amsmath} 
\usepackage{amssymb}
\usepackage{mathrsfs}

\newbox\noforkbox \newdimen\forklinewidth
\setbox0\hbox{$\textstyle\bigcup$}
\setbox1\hbox to \wd0{\hfil\vrule width \forklinewidth depth \dp0
                        height \ht0 \hfil}
\setbox\noforkbox\hbox{\box1\box0\relax}
\def\unionstick{\mathop{\copy\noforkbox}\limits}
  \def\nonfork#1#2_#3{#1\unionstick_{\textstyle #3}#2}
\def\nonforkin#1#2_#3^#4{#1\unionstick_{\textstyle #3}^{\textstyle #4}#2}     
\setbox0\hbox{$\textstyle\bigcup$}
\setbox1\hbox to \wd0{\hfil{\sl /\/}\hfil}
\setbox2\hbox to \wd0{\hfil\vrule height \ht0 depth \dp0 width
                                \forklinewidth\hfil}
\newbox\doesforkbox
\setbox\doesforkbox\hbox{\box1\box0\relax}
\def\nunionstick{\mathop{\copy\doesforkbox}\limits}

\def\fork#1#2_#3{#1\nunionstick_{\textstyle #3}#2}
\def\forkin#1#2_#3^#4{#1\nunionstick_{\textstyle #3}^{\textstyle #4}#2}

\newtheorem{theorem}{Theorem}[section] 
\newtheorem{claim}[theorem]{Claim}

\newtheorem{conclusion}[theorem]{Conclusion}
\newtheorem{observation}[theorem]{Observation}

\theoremstyle{definition}
\newtheorem{definition}[theorem]{Definition}

\newtheorem{fact}[theorem]{Fact}

\newtheorem{conjecture}[theorem]{Conjecture}
\newtheorem{convention}[theorem]{Convention}

\theoremstyle{remark}
\newtheorem{remark}[theorem]{Remark}
\newtheorem{question}[theorem]{Question}
\newtheorem{notation}[theorem]{Notation}

\newtheorem{thesis}[theorem]{Thesis}

\newcommand{\tp}{{\rm tp}}

\newcommand{\ai}{{\rm ai}}

\newcommand{\wed}{{\rm wd}}

\newcommand{\sub}{{\rm sub}}

\newcommand{\Sk}{{\rm Sk}}

\newcommand{\dep}{{\rm dp}}

\newcommand{\odd}{{\rm odd}}

\newcommand{\aec}{{\rm aec}}

\newcommand{\Mod}{{\rm Mod}}
\newcommand{\fp}{{\rm fp}}
\newcommand{\PC}{{\rm PC}}
\newcommand{\MA}{{\rm MA}}

\newcommand{\rt}{{\rm rt}}
\newcommand{\tr}{{\rm tr}}
\newcommand{\ie}{{\rm ie}}
\newcommand{\oor}{{\rm or}}
\newcommand{\sor}{{\rm sor}}

\newcommand{\proj}{{\rm proj}}

\newcommand{\EM}{{\rm EM}}

\newcommand{\pit}{{\rm pit}}
\newcommand{\lev}{{\rm lev}}

\newcommand{\Card}{{\rm Card}}

\newcommand{\lin}{{\rm lin}}
\newcommand{\suc}{{\rm suc}}

\newcommand{\Ord}{{\rm Ord}}

\newcommand{\arity}{{\rm arity}}

\newcommand{\id}{{\rm id}}

\newcommand{\qf}{{\rm qf}}

\newcommand{\ER}{{\rm ER}}
\newcommand{\Dp}{{\rm Dp}}

\newcommand{\LST}{{\rm LST}}

\newcommand{\Dom}{{\rm Dom}}
\newcommand{\Rang}{{\rm Rang}}

\newcommand{\rest}{{\restriction}}

\newcommand{\wilog}{{\rm without loss of generality}}
\newcommand{\Wilog}{{\rm Without loss of generality}}

\newcommand{\then}{{\underline{then}}}
\newcommand{\when}{{\underline{when}}}
\newcommand{\where}{{\underline{where}}}
\newcommand{\Then}{{\underline{Then}}}

\newcommand{\If}{{\underline{if}}}
\newcommand{\Iff}{{\underline{iff}}}
\newcommand{\mn}{{\medskip\noindent}}
\newcommand{\sn}{{\smallskip\noindent}}

\newcommand{\cE}{{\mathscr E}}

\newcommand{\gk}{{\mathfrak k}}

\newcommand{\cC}{{\mathscr C}}

\newcommand{\cI}{{\mathscr I}}

\newcommand{\cL}{{\mathscr L}}
\newcommand{\bbL}{{\mathbb L}}

\newcommand{\bbP}{{\mathbb P}}
\newcommand{\cP}{{\mathscr P}}

\newcommand{\gs}{{\mathfrak s}}

\newcommand{\cS}{{\mathscr S}}

\newcommand{\cT}{{\mathscr T}}

\newcommand{\cX}{{\mathscr X}}
\newcommand{\cY}{{\mathscr Y}}

\newcommand{\cf}{{\rm cf}}
\newcommand{\pr}{{\rm pr}}

\newcount\skewfactor
\def\mathunderaccent#1#2 {\let\theaccent#1\skewfactor#2
\mathpalette\putaccentunder}
\def\putaccentunder#1#2{\oalign{$#1#2$\crcr\hidewidth
\vbox to.2ex{\hbox{$#1\skew\skewfactor\theaccent{}$}\vss}\hidewidth}}

\newenvironment{PROOF}[2][\proofname.]
   {\begin{proof}[#1]\renewcommand{\qedsymbol}{\textsquare$_{\rm #2}$}}
   {\end{proof}}

\begin{document}

\title {A.E.C. with not too many models}
\author {Saharon Shelah}
\address{Einstein Institute of Mathematics\\
Edmond J. Safra Campus, Givat Ram\\
The Hebrew University of Jerusalem\\
Jerusalem, 91904, Israel\\
 and \\
 Department of Mathematics\\
 Hill Center - Busch Campus \\ 
 Rutgers, The State University of New Jersey \\
 110 Frelinghuysen Road \\
 Piscataway, NJ 08854-8019 USA}
\email{shelah@math.huji.ac.il}
\urladdr{http://shelah.logic.at}
\thanks{The author would like to thank the Israel Science Foundation for
partial support of this research (Grant No. 242/03).
The author thanks Alice Leonhardt for the beautiful typing.  Paper 893
in the author's list.  First Version - 2004/June/20.}






\date{December 5, 2013}

\begin{abstract}
Consider an a.e.c. (abstract elementary class), that is, a class $K$ of
models with a partial order refining $\subseteq$ (submodel) which satisfy the
most basic properties of an elementary class.  Our test question is
trying to show that the function $\dot I(\lambda,K)$, counting the
number of models in $K$ of cardinality $\lambda$ up to 
isomorphism, is ``nice", not chaotic, even without assuming it is
sometimes 1, i.e. categorical in some $\lambda$'s.  We prove
here that for some closed unbounded class $\bold C$ of
cardinals we have (a),(b) or (c) where
\mn
\begin{enumerate}
\item[$(a)$]  for every $\lambda \in \bold C$ of cofinality $\aleph_0,
\dot I(\lambda,K) \ge \lambda$ 
\sn
\item[$(b)$]  for every $\lambda \in \bold C$ of cofinality $\aleph_0$ and 
$M \in K_\lambda$, for every cardinal $\kappa \ge \lambda$ there is
$N_\kappa$ of cardinality $\kappa$ extending $M$ (in the sense of our
a.e.c.)
\sn
\item[$(c)$]  $\gk$ is bounded; that is, $\dot I(\lambda,K) =0$ for
  every $\lambda$ large enough (equivalently $\lambda \ge
  \beth_{\delta_*}$ where $\delta_* = (2^{\LST(\gk)})^+$).
\end{enumerate}
\mn
Recall that an important difference of non-elementary classes from the
elementary case is the possibility of having models in $K$, even of
large cardinality, which are maximal, or just failing clause (b).
\end{abstract}

\maketitle
\numberwithin{equation}{section}
\setcounter{section}{-1}
\newpage

\centerline {Anotated Content}
\bigskip

\noindent
\S0 \quad Introduction to the subject, pg.\pageref{Intro}
\bigskip

\noindent
\S1 \quad Introduction to the paper, pg.\pageref{Introduction}

\S(1A) \quad Content, pg.\pageref{content}

\S(1B) \quad Discussion, (label y), pg.\pageref{Discussion}

\S(1C) \quad What is done, pg.\pageref{What}

\S(1D) \quad Recalling Definitions and Notation, (label z), pg.\pageref{Defno}
\mn
\begin{enumerate}
\item[${{}}$]   [Defining $\Upsilon^{\oor}_\kappa[\gk],\Upsilon^{\sor}$.]
\end{enumerate}
\bigskip

\noindent
\S2 \quad More on Templates, pg.\pageref{More}
\mn
\begin{enumerate}
\item[${{}}$]  On $\gk_\mu,\Upsilon^{\sor}_\kappa[\gk]$, the ideal ER,
being standard; the ideal ER, \cite{Sh:318}; isomorphic of
vocabularies, \ref{z19} - \ref{z21}; partial orders on
$\Upsilon^{\oor}_\kappa[\gk]$, \ref{z21}, \ref{z24}, basic properties
of $\le^3_\kappa,\le^4_4$, see \ref{z26}.  On
$\Upsilon^{\sor}_\kappa[k_M] \ne \emptyset$, \ref{z29}; amalgamating
$\le^\oplus_\kappa,\le^4_\kappa$.  Lastly, we define $\pit(\cT,\bold
I)$, \ref{z32} and have the relevant partition theorem, \ref{z35}.
\end{enumerate}
\bigskip

\noindent
\S3 Approximations to EM models, (label a) , pg.\pageref{Approx}
\mn
\begin{enumerate}
\item[${{}}$]   We define direct/pre-witnesses, (\ref{a5}, \ref{a6}),
prove existence, \ref{a9}.  Deal with indirect witnesses (used in the
main proof (see \ref{a12}, \ref{a14}) and prove the main result.
\end{enumerate}
\bigskip

\noindent
\S4 \quad Concluding Remarks, pg.\pageref{Concluding}
\newpage

\section {Introduction to the subject} \label{Intro}

We would like to have classification theory for non-elementary classes
$K$ and more specifically to generalize stability.  Naturally we 
use the function $\dot I(\lambda,K) =$ 
number of models up to isomorphism, as a major test
problem.  Now ``non-elementary" has more than one interpretation, we
shall start with the infinitary logics $\bbL_{\lambda,\kappa}$.

There are other directions; mostly where compactness in some form
holds (e.g. a.e.c. with amalgmation, see about those in \cite{Sh:E53},
and on a try to blend with descriptive set theory see \cite{Sh:849}).  We had
held that for $\kappa > \aleph_0$ the above cannot be developed as, e.g. if
$\bold V = \bold L$ or just $\bold V \models ``0^\#$ does not exist",
then there is $\psi \in \bbL_{\aleph_1,\aleph_1}$ such that if
$\cf(\mu) = \aleph_0 \wedge (\forall \alpha < \mu)(|\alpha|^{\aleph_0}
< \mu)$ \then \, $M \models \psi,\|M\| = \mu$ \Iff \, $M \cong
(\bold L_\mu,\in)$.  However, lately \cite{Sh:F1228} gives evidence that
for $\theta$ a compact cardinal, we can
generalize to $\bbL_{\theta,\theta}$ some theorems of
\cite[Ch.VI]{Sh:c} on saturation of ultra-powers and Keisler's order.
This shows that stability theory for $T \subseteq
\bbL_{\theta,\theta}$ exists, but it is still not clear how far we can go
including $A = |N|,N \prec M$ and $\cup\{M_u:u \subset n\}$ when
$\langle M_u:u \subset n\rangle$ is a so called stable $\cP^-(n)$-system.

Anyhow (for the purposes of this history, and the present paper) 
we now concentrate on $\Mod_\psi,\psi \in
\bbL_{\lambda^+,\aleph_0}$ so $\kappa =\aleph_0$.  Here we have both
downward $\LST$ theorems, even using $\le \lambda$ finitary Skolem
functions.  Also we have the upward $\LST$ theorem, using $\EM$
models.

Naturally all works started with assuming categoricity in some cardinal,
except some dealing with the $\aleph_n$'s for $\psi \in
\bbL_{\aleph_1,\aleph_0}$.  In this case we may many times deal with
$\psi \in \bbL_{\aleph_1,\aleph_0}(Q)$.  Some works apeared in the
eighties (see the books \cite{Bal09}, also \cite{Sh:h}, \cite{Sh:i}).  

\begin{definition}
\label{y9}
Let $\dot I(\lambda,K)$ be the cardinality of 
$\{M /\cong:M \in K$ of cardinality
$\lambda\}$ where $K$ is a class of $\tau(K)$-models (e.g. $K = K_{\gk}$
where $\gk = (K_{\gk},\le_{\gk})$).
\end{definition}

First, in ZFC, answering a question of Baldwin, it was proved that
 $\psi$ cannot be categorical,
moreover if $\dot I(\aleph_1,\psi)=1$ then $\dot I(\aleph_2,\psi) \ge 1$.
Also if $\dot I(\aleph_1,\psi) < 2^{\aleph_1}$, then for some
countable first order $T$ with an atomic model $K_T = \{M:M$ an atomic
model of $T\}$ is $\subseteq \Mod_\psi$, but $1 \le \dot I(\aleph_1,K_T)$.
Fix $T$ for awhile, now if $2^{\aleph_n} < \aleph_{n+1},\dot I(\aleph_n,T) <
\mu_{\wed}(\aleph_{n+1},2^{\aleph_n})$ for\footnote{note that
  $\mu_{\wed}(\lambda^+,2^\lambda)$ is essentially $2^{\lambda^+}$.}
every $n$ then $K_T$ is excellent which means it is quite similar to
the class of models of an $\aleph_0$-stable countable complete first
order theory.  For this we consider $\bold S^m(A,M)$ for $A \subseteq M \in
K_T$, only some ``nice" $A$.  On the other hand for any $n$ for some such
$T_n,K_{T_n}$ is categorical in every $\lambda \le \aleph_n$ but $\dot
I(\lambda,T) = 2^\lambda$ for $\lambda$ large enough.  However, we do
not know:

\begin{conjecture}
\label{x4}
(Baldwin)
If $K_T$ is categorical in $\aleph_1$, \then \, $K_T$ is
$\aleph_0$-stable, equivalently is absolutely categorical.
\end{conjecture}

\noindent
Related is the:
\begin{conjecture}
\label{x6}
If $K_T$ is categorical in $\aleph_1$ but not $\aleph_0$-stable \then
\, $\dot I(2^{\aleph_0},K_T) = \beth_2$.

See work in preparation Baldwin-Laskowski-Shelah (\cite{Sh:F1098}) on
such $K_T$'s; it certainly says there is a positive theory for such
classes (e.g. pseudo minimal types exist).  We recently have changed
our mind and now think:
\end{conjecture}

\begin{conjecture}
\label{x9}
If $K_T$ is categorical in every $\aleph_n$ \then \, $K_T$ is
excellent.

This means that the present counter-examples are best possible.  As
this seems very far we may consider a weaker conjecture.
\end{conjecture}

\begin{conjecture}
\label{x12}
Assume $\bbP$ is a c.c.c. forcing notion of cardinality $\lambda$ such
that $\Vdash_{\bbP} ``\MA + 2^{\aleph_0} = \lambda"$ and $\lambda =
\lambda^{< \lambda} > \beth_{\omega_1}$.  
If $K_T$ is categorical in every $\lambda <
2^{\aleph_0}$ then $K_T$ is excellent.

There is more to be said, see \cite{Sh:F1273}. 
\end{conjecture}
\bigskip

\centerline {$* \qquad * \qquad *$}
\bigskip

In another direction, the investigation of models of cardinality
$\aleph_1$ does not point to a canonical choice of logic for which
the theorems on $\dot I(\psi,\aleph_1) = 1$ 
holds.  This had motivated the definition of a.e.c. $\gk =
(K_{\gk},\le_{\gk})$ which has the ``bottom" property of elementary
class $K = (\Mod_T,\prec),T$ a complete first order theory
(i.e. $K_{\gk}$, a class of $\tau_{\gk}$-models, $\le_{\gk}$ a partial
order on it, both closed under isomorphism, union under
$\le_{\gk}$-directed systems of member of $K_{\gk}$ belong to
$K_{\gk}$, moreover is a $\le_{\gk}$-lub (= union of a directed system
of $\le_{\gk}$-submodels of $N$ is a $\le_{\gk}$-submodel of $N$),
existence of a LST number and $M_1 \subseteq M_2 \wedge M_1 \le_{\gk}
N \wedge M_2 \le_{\gk} N \Rightarrow M_1 \le_{\gk} M_2$).
\bigskip

\begin{thesis}
\label{x12}
1) The framework of a.e.c. $\gk$ is wider and not too far and better than
   the family of $(\Mod_\psi,\prec_{\sub(\psi)})$ where $\psi \in
   \bbL_{\lambda^+,\aleph_0}$. 

\noindent
2) The right generalization of types in this context is orbital types.  

Why?  The ``wider" in \ref{x12}(1) is obvious.  The ``not
   too far" is by the representation theorem which says that for some
   vocabulary $\tau_1 \supseteq \tau(\gk)$ of cardinality 
$\le \lambda,\lambda$ the $\LST$-number $+|\tau(\gk)|$ 
and set $\Gamma$ of quantifier free 1-types, 
$K_{\gk} = \PC(\emptyset,\Gamma) = \{M \rest \tau_{\gk}:M$
   a $\tau_1$-model omitting every $p(x) \in \Gamma\}$; similarly
   $\le_{\gk}$.  We can deduce the upward $\LST$, and so existence of
   suitable $\Phi \in \Upsilon^{\lin}[\gk]$ so we have $\EM$-models.
   For $\gk$ with $\LST_{\gk} = \aleph_0$ it is natural to restrict
   ourselves to the case ``$\Gamma$ is countable" above for both
   $K_{\gk} \le_{\gk}$, then we say $\gk$ is $\aleph_0$-presentable.
   So we may wonder for such $\gk$ \If \, $n< \omega \Rightarrow
   2^{\aleph_n} + \dot I(\aleph_{n+1},K_{\gk}) <
   \mu_{\wed}(\aleph_{n+1},2^{\aleph_n})$ implies $\gk$ satisfies the
   parallel of being excellent?  The answer is yes by \cite{Sh:h},
   \cite{Sh:i}, but the way is long.  Also, we may replace $\aleph_0$
   by any $\lambda$ provided that $I(\lambda,K_{\gk}) = 1 =
   I(\lambda^+,K_{\gk})$ and $1 \le \dot I(\lambda^{++},K_{\gk}) <
   \mu_{\wed}(\lambda^{++},2^{\lambda^+})$, see more in \cite{Sh:E53}.

A central notion there is ``$\gs$ is a good $\lambda$-frame",
$\gk_{\gs} = \gk,\LST_{\gk} \le \lambda$, this is ``bare bones
superstable".

This is enough for proving
\mn
\begin{enumerate}
\item[$(*)$]  if ($\gk$ is an a.e.c.), $\LST_{\gk} \le
  \lambda,2^{\lambda^{+n}} < 2^{\lambda^{+n^+}}$ and
  $\dot I(\lambda^{+n},K_{\gk})=1$ for every $n$ and $K_{\gk}$ has
  models of cardinality $\ge \beth_{(2^{\LST(\gk)})^+}$, \then \,
  $K_{\gk}$ is categorical in every $\mu \ge \lambda$.
\end{enumerate}
\end{thesis}

\noindent
However
\begin{conjecture}
\label{x15}
If $\gk$ is an a.e.c., $K_{\gk}$ is categorical in some $\lambda$
large enough than $\LST_{\gk}$, \then\, $K_{\gk}$ is categorical in
every $\mu \ge \lambda$.

Note that \cite{Sh:734} is a step ahead: in the context of \ref{x15},
for many $\mu = \beth_\mu \in [\LST_{\gk},\lambda)$, there is a good
  $\mu$-frame $\gs_\mu$ such that $\gk_{\gs} = K^{\gk}_\mu$.  If we
  have this for $\omega$ successive $\mu$'s we shall be done by
  \cite{Sh:600}, but in \cite{Sh:734} the family of such $\mu$'s is
scattered; a beginning is \cite{Sh:842}.
\end{conjecture}

\noindent
A much harder conjecture is:
\begin{conjecture}
\label{x17}
1) The main gap theorem holds for a.e.c. $K_{\gk}$ for $\lambda$ large
   enough.

\noindent
2) The class $\sup-\lim_{\gk} = \{\lambda$: there is a super-limit $M
   \in K^{\gk}_\lambda\}$ is ``nice", e.g. contains every large enough
   $\lambda$ \underline{or} contains no large enough $\lambda$.
\end{conjecture}
\bigskip

\centerline {$* \qquad * \qquad *$}
\bigskip

We may wonder
\begin{question}
\label{x23}
1) Maybe there is a natural logic which is the natural framework for
   categoricity spectrum.

\noindent
2) Also for the super-limit spectrum.

We expect such logic to be stronger than $\bbL_{\lambda^+,\aleph_0}$
but weaker than $\bbL_{\lambda,\lambda}$.  This may remind us of
\cite{Sh:797}.  The logic discovered there is $\bbL^1_{<\lambda}$ for
$\lambda = \beth_\lambda$, it is between $\bbL^{-1}_{<\lambda} =
\cup\{\bbL_{\mu^+,\aleph_0}:\mu < \lambda\}$ and $L^0_{<\lambda,\mu} =
\cup\{\bbL_{\mu^+,\mu^+}:\mu < \lambda\}$, in a strong way well
ordering is not well defined and it can be characterized (as
Lindstr\"om theorem characterize first order logic) and has
interpolation.   In addition,
for $\lambda$ a compact cardinal $\bbL^1_{< \lambda}$-equivalence of
$M_1,M_2$ is equivalent to having isomorphism $\omega$-limit
ultra-powers by $\lambda$-complete ultrafilters, see \cite{Sh:F1228}.

However, probably the characterization in \cite{Sh:797} was by ``the
maximal logic such that ...".  So maybe we should restrict the logic further
such that ``$\EM$ model can be constructed".

We conjecture there is a logic characterized by being maximal under
this stronger demand, and in it we can say at least something on the
function $\dot I(\lambda,\psi)$, and maybe much.  This is interesting
also from the point of view of soft model theory: we conjecture that
there are many such intermediate logics with characterization (and the
related interpolation theorem).
\end{question}
\newpage

\section {Introduction to the paper} \label{Introduction}

In this section, we begin by motivating our line of investigation.  
See notation in \S(1D) below 
(and more self contained introduction in \S(1B), \S(1C)).
\bigskip

\subsection {Motivation/Content} \label{content}\
\bigskip

\noindent
We knew of old (see: \cite[Ch.XIII,4.15]{Sh:c}):
\begin{theorem}
\label{y1} 
For a countable complete \underline{first order} theory $T$,
one of the following holds:
\mn
\begin{enumerate}
\item[$(a)$]  $T$ is categorical in every $\lambda > \aleph_0$
\sn
\item[$(b)$]  $\dot I(\lambda,T) = \beth_2$ for every cardinal $\lambda \ge
  2^{\aleph_0}$
\sn
\item[$(c)$]  $\dot I(\aleph_\alpha,T) \ge 1 + |\alpha|$ for every
  ordinal $\alpha$.
\end{enumerate}
\end{theorem}

\noindent
For a.e.c. we have something when $\gk$ is categorical in some
$\lambda$'s (\cite{Sh:734}, \cite{Sh:600}) and something about $\dot
I(\aleph_1,\gk)$, (\cite{Sh:88r}, about when $1 \le \dot I(\aleph_1,\gk)
< 2^{\aleph_1}$, particularly when $2^{\aleph_0} < 2^{\aleph_1}$ and
then on higher cardinals) but nothing for general a.e.c. $\gk$.  
The current paper is motivated by hopes of finding something like \ref{y1} for
a.e.c.'s.  Recall the history.

\noindent
Our approach here assumes/relies on:
\begin{thesis}  
\label{y2}
Reasonable to concentrate on cardinals from $\bold C_{\fp} =
\{\lambda:\lambda = \beth_\lambda\}$, where $\fp$ stands for ``fixed points".

Why?  If $\lambda \in \bold C_{\fp},\lambda > \LST(\gk)$ and $M \in
K^{\gk}_\lambda$ \then \, for every $\theta \in [\LST(\gk),\lambda)$
  and $N \le_{\gk} M,\|N\| = \theta$ there is $\Phi \in
  \Upsilon_{\gk,\theta}$ so $|\tau(\Phi)| = \theta$ such that for any linear
  order $I$, e.g. $I = \lambda$ we have 
$N \le_\gk \EM_{\tau(\gk)}(I,\Phi)$.  So
in $K^{\gk}_\lambda$ we have many models of the form
  $\EM_{\tau(\gk)}(I,\Phi),\Phi \in \Upsilon_{\gk,< \lambda}$.  If
  $\dot I(\lambda,\gk) < \lambda$, many of them will be isomorphic.
Hence for many $\theta_1 < \theta_2 < \lambda,\theta_1 \ge \LST(\gk)$,
every $N \le_{\gk} M$ of cardinality $\theta_2$ can be
$\le_{\gk}$-embedded into some $\EM_{\tau(\gk)}(\lambda,\Phi),\Phi \in
\Upsilon^{\oor}_\kappa[\gk]$. 

Informally, the point is it allows us to use EM models.  
The key point is finding a suitable template,
set $\Phi$ of quantifier free types, which requires finding enough
indiscernible sequences.  When $K_{\gk}$ is an a.e.c. (as opposed to an
elementary or pseudo elementary class) we must go through the
Presentation Theorem to find an indiscernible sequence, i.e. we
require sufficiently large models omitting the types in $\Gamma$.

To further motivate our approach, 
consider a not so strong conjecture, still enough to exemplify
``the function $\lambda \mapsto \dot I(\lambda,\gk)$ cannot be too wild".
\end{thesis}

\begin{conjecture}
\label{y4}
1) Letting $\bold C^{\fp}_{\aleph_0} = \{\lambda:\lambda =
\beth_\lambda$ and $\cf(\lambda) = \aleph_0\}$ and fixing an a.e.c. $\gk$, not
both of the following classes are stationary (or restrict yourself to
some strongly inaccessible $\mu$ and ``stationary" means below it):
\mn
\begin{enumerate}
\item[$(a)$]  $\bold S_1 = \{\lambda \in \bold C^{\fp}_{\aleph_0}:\dot
  I(\lambda,\gk) < \lambda\}$
\sn
\item[$(b)$]    $\bold S_2 = \{\lambda \in \bold C^{\fp}_{\aleph_0}:
\dot I(\lambda,\gk) \ge \lambda\}$.
\end{enumerate}
\mn
2) A weaker conjecture (presented in the abstract) is replacing clause
(b) by (h) $\bold S_3 = \{\lambda \in \bold C^{\fp}_{\aleph_0}$: for every
$M \in K^{\gk}_\lambda$ has $\le_{\gk}$-extensions $N$ of any
cardinality $> \lambda\}$.
 
Why ``$\cf(\lambda) = \aleph_0$"?   First, trying to prove $\lambda
\in \bold S_3$, we can approximate $N$ by $\Phi \in
\Upsilon^{\oor}_{\lambda_n}[\gk],\lambda_n < \lambda$ as we can
approximate $M$ by $N' \le_{\gk} M,\|N'\| = \lambda_n$ where
$\lambda_n < \lambda_{n+1} < \lambda = \Sigma\{\lambda_m:m\}$.
Second, for $\lambda \in \bold C^{\fp}_{\aleph_0}$ it is enough to show that
$\{M/\equiv_{\bbL_{\infty,\lambda}}:M \in K^{\gk}_\lambda\}$ is small
because it is well known that if $\cf(\lambda) = \aleph_0$ and
$M_1,M_2$ are of cardinality $\lambda$ and
$\bbL_{\infty,\lambda}$-equivalent then they are isomorphic; on such
logics see, e.g. \cite{Dic85}.
\end{conjecture}

\begin{thesis}
\label{y12}
There are, for a.e.c. $\gk$, meaningful dichotomy theorems for $\dot
I(\lambda,K_{\gk})$ when $K$ is a class of $\tau(\gk)$-models, $K =
K_{\gk}$ and $\gk = (K_{\gk},\le_{\gk})$.

This is a more concrete thesis than ``considering a.e.c.'s is a good frame
for model theory"; even more concrete is the ``main gap conjecture".  It
had been proved that if $K_{\gk}$ is the class of models of a complete
countable first order theory \then \, it satisfies the ``main gap", i.e. either
$\dot I(\lambda,K)$ is large, even 
$= 2^\lambda$ for all uncountable $\lambda$ or $\dot
I(\aleph_\alpha,K)$ is small, even 
$< \beth_{\omega_1}(|\alpha|)$ for all $\alpha > 0$; 
see \cite[Ch.XII]{Sh:c}, ``The book's main theorem".  In general for a
class $K$ of $\tau$-models the ``main gap" will say that either 
$\dot I(\lambda,K)$ is large (i.e. $2^\lambda$ or $\ge \lambda^+$) 
for every $\lambda$ large enough or it is small for every $\lambda$ large
enough say $\dot I(\aleph_\alpha,K)$ is $\le \beth_{1,n}(|\alpha|)$
for some $n=n(K) < \omega$.

We are far away from this, still, until now for the a.e.c. the
categoricity case was almost alone, i.e. we start assuming
$\dot I(\lambda,K)=1$ in some $\lambda$, see below,
but we try here to look ``higher".  

The contribution of the present
paper is to show that in the much more general context of a.e.c.'s for
some $\aleph_0$-closed unbounded class $\bold C$ of cardinals, we have
$\lambda \in \bold C \Rightarrow \dot I(\lambda,K_{\gk}) \ge \lambda$,
a non-structure result, \underline{or} $\lambda \in \bold C \wedge M
\in K^{\gk}_\lambda \Rightarrow M$ has arbitrary large
$\le_{\gk}$-extensions.  Note that the latter property is now taken for
granted for elementary classes but is a real gain for a.e.c.

As noted in \S0, in \cite{Sh:734} and
\cite{Sh:600} we obtained results on $\dot I(\lambda,K)$ for a.e.c.'s
assuming categoricity in some $\lambda$'s.  
However, nothing was known for general a.e.c.'s under
weaker few models assumption.

On abstract elementary classes, see \cite{Sh:88r}, \cite{Bal09} and
 \cite{Sh:E53}.  We will make essential use of
the Presentation Theorem, which says that every a.e.c. can be
represented as a PC class, say $\PC(T,\Gamma)$, see
 \cite[\S1]{Sh:88r}.

We thank the audience in the lecture in the Hebrew University seminar
2/2005 for their comments on an earlier version of this paper and
Maryanthe Malliaris for helping much in improving \S1 and some
corrections in fall 2011 - winter 2012 and Will Boney for some further
corrections (fall 2013).
\bigskip
\end{thesis}
\bigskip

\subsection {Discussion} \label{Discussion}\
\bigskip

We give some further details regarding \S(1A).

Why?  In Thesis \ref{y2} the result on EM models needed
is: \cite[Claim 0.6]{Sh:394}, \cite[Claim 8.6]{Sh:394}, 
the ``a.e.c. omitting types theorem" and
\cite[Lemma 8.7,p.46]{Sh:394}.

\begin{fact}
\label{y15}
Let $\gk$ be an a.e.c.  If $\lambda \in \bold C_{\fp},\lambda >
\LST_{\gk}$ and $M \in K^{\gk}_\lambda$ \then \, 
for every $\theta \in [\LST_{\gk},\lambda)$ and 
$N \le_{\gk} M$ of cardinality $\theta$ there is $\Phi \in
\Upsilon[\gk]$ such that:
\mn
\begin{enumerate}
\item[$(a)$]  $|\tau(\Phi)| =\theta$
\sn
\item[$(b)$]  for any linear order $I$, in particular $I = \lambda$,
\wilog \, $N \le_{\gk} \EM_{\tau(\gk)}(I,\Phi)$ where this
  denotes the reduct of the EM model to the vocabulary of $\gk$. 
\end{enumerate}
\bigskip

\noindent
\underline{Comment}:

Let us repeat, the two points when $\cf(\lambda) = \aleph_0$ may be as
required:
\mn
\begin{enumerate}
\item[$(a)$]  downward large depth in \S3,
\sn
\item[$(b)$]  if we like to find large $N \le_{\gk}$-extending $M$ 
for a given $M \in
  K^{\gk}_\lambda$, if $\cf(\lambda) = \aleph_0$ we can get it as an 
$\omega$-limit of $M' <_{\gk} M,\|M'\| < \lambda$.
\end{enumerate}
\end{fact}

\noindent
Such considerations further lead us to
\begin{question}
\label{y6}
Let $\Phi \in \Upsilon_\theta[\gk]$ and $\kappa$ be a cardinal.

Sort out the functions
\mn
\begin{enumerate}
\item[$(a)$]   $\lambda \mapsto |\{\EM_{\tau(\gk)}(I,\Phi)/\cong:I$ a
linear order of cardinality $\lambda\}|$
\sn
\item[$(b)$]  $\lambda \mapsto \dot I_{\tau(\gk)}(\lambda,\kappa,\Phi) :=
  |\{\EM_{\tau(\gk)}(I,\Phi)/\equiv_{\bbL_{\infty,\kappa}}:I$ a linear
  order of cardinality $\lambda\}$.
\end{enumerate}
\mn
Recall, by \cite{Sh:11} restricting ourselves to cardinals 
$\lambda = \lambda^{< \kappa}$, that the function in clause (b) of \ref{y6} 
is ``nice", more specifically: if 
$\theta \le \lambda_1 = \lambda^{< \kappa}_1 <
\lambda_2$ \then \, $\dot I_{\tau(\gk)}(\lambda_1,\kappa,\gk) \ge
  \min\{\lambda^+_1,\dot I(\lambda_2,\kappa,\gk)\}$.

What occurs if $\lambda_1 < \lambda^{< \kappa}_1$?  The case $\lambda_1 =
\beth_\delta,\cf(\delta) = \aleph_0$ is more approachable than the
general case, see \ref{n4}.

Our hope is to get ``bare bones superstability", i.e. good $\lambda$-frames
inside $\gk$, (as in \cite{Sh:600},\cite{Sh:734}).

Another point concerning the function 
$\dot I(\lambda,\kappa,\gk)$ is: for a model
$M$, cardinal $\theta$ and logic $\cL$ we can define the depth of $M$
for $(\cL,\theta)$ as $\min\{\alpha$: if $\bar a,\bar b \in
{}^\varepsilon M,\varepsilon < \theta$ and $\bar a,\bar b$ realizes the
same formulas of $\bbL_{\infty,\theta}$ (or
$\bbL_{\infty,\theta}[\gk])$ of depth $< \alpha$ then they realize the
same $\bbL_{\infty,\theta}$-formulas$\}$; of course, only formulas in
$L_{\|M\|^{< \theta},\theta}$ are relevant.  This is a good way to
``slice" the equivalence and it is easier for LST considerations.
\end{question}
\bigskip

\subsection {What is Done} \label{What}\
\bigskip

A phenomena making the investigation of general
a.e.c. hard is having $\le_{\gk}$-maximal
models of large cardinality.  As with amalgamation, we may consider the
property
\mn
\begin{enumerate}
\item[$(*)^1_\lambda$]  if $M \in K^{\gk}_\lambda$ then $M$ is not
  $\le_{\gk}$-maximal.
\end{enumerate}
\mn
In investigations like \cite{Sh:E46} and \cite{Sh:576}, which look at
$\cup\{K^{\gk}_{\lambda^{+ \ell}}:\ell < 4\}$ this is relevant.  But
in investigations as in \cite{Sh:734}, looking at
$\cup\{K^{\gk}_\lambda:\lambda = \beth_\lambda\}$, it is more natural
to consider
\mn
\begin{enumerate}
\item[$(*)^2_\lambda$]  if $M \in K^{\gk}_\lambda$ then for any $\mu >
  \lambda$ there is $N \in K^{\gk}_\mu$ which $\le_{\gk}$-extends $N$.
\end{enumerate}
\mn
In \S3 we consider a $\lambda = \beth_\lambda$ 
of cofinality $\aleph_0$ which is
more than strong limit and try to prove non-structure from
$\neg(*)^2_\lambda$. Given $N \in K^{\gk}_\lambda$ we try to build an
EM model (that is construct the $\Phi$) $\le_{\gk}$-extending $N$ by
an increasing chain of approximations: given $\lambda_n \rightarrow
\lambda,M_n \rightarrow N,M_n
\in K^{\gk}_{\lambda_n}$. The $n$-th approximation $\Phi_n$ to
$\Phi$ has to have ``$\Phi_n$ in a suitable sense is
represented in $N$ say of size $\lambda_{n+1}"$.

Being stuck should be a reason for non-structure.  For simplicity we
consider only cardinals $\mu = \beth_\mu$, the gain
without this restriction seems minor.

Concerning the results of \S3 it would be nicer to make one more step
concerning \ref{a25}, \ref{a22} and deal also with $\lambda =
\beth_\lambda$ instead of $\lambda = \beth_{1,\lambda}$, 
but a more central question is to get the
non-structure result for every $\lambda' > \lambda$.  It is natural to
try given $\Phi \in \Upsilon^{\sor}_\kappa[\gk_M]$ and $M \le_{\gk}
N$, to define a
``depth" for approximation of the existence of a $\le_{\gk}$-embedding
of standard $\EM_{\tau(\gk)}(I,\Phi)$ into $N$ (see Definition
\ref{z14}(2)), so that depth infinity give existence.  But this does not
work for us, so Definition \ref{a5} is a substitute, moreover we need
``indirect evidence", see Definition \ref{a12}.

\noindent
Our main theorem is
\begin{theorem}
\label{y23}
For any a.e.c. for some closed unbounded class of cardinals $\bold C$,
if $(\exists \lambda \in \bold C)[\cf(\lambda) = \aleph_0 \wedge \dot
I(\lambda,K_{\gk}) < \lambda]$ and $M \in K_{\gk}$ of cardinality $\mu
\in \bold C$ of cofinality $\aleph_0,M$ has a proper
$<_{\gk}$-extension, and even ones of arbitrarily large cardinality.
\end{theorem}

\noindent
The natural next steps are
\begin{conjecture}
\label{y10}
1) In Theorem \ref{a28}, i.e. what is promised in the abstract we can
   choose $\bold C$ as an end segment of $\{\mu:\mu = \beth_{1,\mu}\}$
   or just choose $\bold C$ as $\{\mu:\mu = \beth_{2,\mu}\}$.

\noindent
2) For every a.e.c. $\gk$ for some closed unbounded class $\bold C$ of
cardinals, we have $M \in K^{\gk}_\lambda \wedge \lambda \in \bold C
\wedge \cf(\lambda) = \aleph_0 \Rightarrow 
\Upsilon^{\oor}_\lambda[\gk_M] \ne \emptyset$ \underline{or}
   $\lambda \in \bold C \wedge \cf(\lambda) = \aleph_0 \Rightarrow
   \dot I(\lambda,K_{\gk}) \ge 2^\lambda$ or at least $\ge \lambda^+$.

We intend to deal with part (1) in a continuation.
\end{conjecture}
\bigskip

\subsection {Recalling Definitions and Notation} \label{Defno} \
\bigskip

\begin{notation}
\label{z0}
Let $\Card$ be the class of infinite cardinals.
\end{notation}

\begin{definition}
\label{z1}
1) Let $\beth_{0,\alpha}(\lambda) =
\beth_\alpha(\lambda) := \lambda + \Sigma\{2^{\beth_\beta(\lambda)}:\beta 
< \alpha\}$.  Let $\beth_{\varepsilon,\alpha}(\lambda)$ 
be defined by induction on $\varepsilon > 0$ and for each $\varepsilon$ by
induction on $\alpha:\beth_{\varepsilon,0}(\lambda) = 
\lambda$, for limit $\beta$ we let
$\beth_{\varepsilon,\beta} = \sum\limits_{\gamma < \beta}
\beth_{\varepsilon,\gamma}$ and for $\varepsilon = \zeta +1$ let
$\beth_{\zeta +1,\beta + 1}(\lambda) = \beth_{\zeta,\mu}$
where\footnote{why not, e.g. $\mu = \beth_{1,\beta}(\lambda)^+$?  Not
a serious difference as for limit $\alpha$ we shall get the same value
and in \ref{z8}(1) this simplifies the notation.}
$\mu = (2^{\beth_{\zeta,\beta}(\lambda)})^+$, lastly for limit
$\varepsilon$ let $\langle \beth_{\varepsilon,\alpha}:\alpha \in
\Ord\rangle$ list in increasing order the closed unbounded class
$\bigcap\limits_{\zeta < \varepsilon}\{\beth_{\zeta,\alpha}:\alpha \in \Ord\}$.

\noindent
2) Let $\lambda \gg \kappa$ mean $(\forall \alpha <
   \lambda)(|\alpha|^\kappa < \lambda)$.
\end{definition}

\begin{convention}
\label{z2}
1) $\gk = (K_{\gk},\le_{\gk})$ is an a.e.c., with vocabulary
$\tau_{\gk} = \tau(\gk)$ and $\LST(\gk) = \LST_{\gk}$ 
its L\"owenheim-Skolem-Tarski number,
see \cite[\S1]{Sh:88r}.  If not said otherwise, we assume
$|\tau_{\gk}| \le \LST_{\gk}$.

\noindent
2) $K^{\gk}_\lambda = K_{\gk,\lambda} = \{M \in K_{\gk}:\|M\| =
   \lambda\}$.

\noindent
3) If $K = K_{\gk}$ we may write $\gk$ instead of $K$; also we may
 write $K$ or $K_\lambda$ omitting $\gk$ when (as usually here) 
$\gk$ is clear from the context.
\end{convention}

\begin{definition}
\label{z3}
For a class $K$ of $\tau$-models:
\mn
\begin{enumerate}
\item[$(a)$]  for a cardinal $\lambda$, let 
$\dot I(\lambda,K)$ be the cardinality of 
$\{M/\cong:M \in K$ has cardinality $\lambda\}$
\sn
\item[$(b)$]  for a cardinal $\lambda$ and a logic $\cL$, let $\dot
  I(\lambda,\cL,K) = \{M/\equiv_{\cL(\tau)}:M \in K$ has cardinality
  $\lambda\}$.
\end{enumerate}
\end{definition}

\begin{definition}
\label{z4}
1) $\Phi$ is a template proper for linear orders \when \,:
\mn
\begin{enumerate}
\item[$(a)$]   for some vocabulary $\tau = \tau_\Phi =
\tau(\Phi),\Phi$ is an $\omega$-sequence, with the $n$-th element a complete
quantifier free $n$-type in the vocabulary $\tau$,
\sn
\item[$(b)$]   for every linear order $I$ there is a $\tau$-model
$M$ denoted by $\EM(I,\Phi)$, generated by $\{a_t:t \in I\}$ such that
$s \ne t \Rightarrow a_s \ne a_t$ for $s,t \in I$ and
 $\langle a_{t_0},\dotsc,a_{t_{n-1}}\rangle$ realizes the quantifier free
$n$-type from clause (a) whenever $n < \omega$ and
$t_0 <_I \ldots <_I t_{n-1}$.  We call $(M,\langle a_t:t \in
I\rangle)$ a $\Phi-\EM$-pair or $\EM$-pair for $\Phi$; 
so really $M$ and even $(M,\langle a_t:t \in I\rangle)$ are  
determined only up to isomorphism but abusing notation we may ignore this and
use $I_1 \subseteq J_1 \Rightarrow \EM(I_1,\Phi) \subseteq 
\EM(I_2,\Phi)$.  We call $\langle a_t:t \in I\rangle$ ``the"
skeleton of $M$; of course again ``the" is an abuse of notation as it is not
necessarily unique.
\end{enumerate}
\mn
1A) If $\tau \subseteq \tau(\Phi)$ then we let $\EM_\tau(I,\Phi)$ be the
$\tau$-reduct of $\EM(I,\Phi)$.

\noindent
2) $\Upsilon^{\oor}_\kappa[{\gk}]$ is the class of templates $\Phi$
proper for linear orders satisfying clauses 
$(a)(\alpha),(b),(c)$ of Claim \ref{z8}(1) below and
$|\tau(\Phi) \backslash \tau_{\gk}| \le \kappa$; normally we assume
$\kappa \ge |\tau_{\gk}| + \LST_{\gk}$ but using $\gk_M$ we do not
assume $\kappa \ge \|M\|$, see \ref{z12}.  The default value of $\kappa$ is
$\LST_{\gk}$ and then we may write $\Upsilon^{\oor}_{\gk}$
or $\Upsilon^{\oor}[{\gk}]$ and for simplicity if not said otherwise
$\kappa \ge \LST_{\gk}$ (and so $\kappa \ge |\tau_{\gk}|$).  We may omit $\gk$
when clear from the context.

\noindent
3) For a class $K$ of so called index models, we 
define ``$\Phi$ proper for $K$" similarly when in clause (b) of
part (1) we demand $I \in K$, so $K$ is a class of $\tau_K$-models,
i.e.
\mn
\begin{enumerate}
\item[$(a)$]   $\Phi$ is a function, giving for any complete quantifier free
$n$-type in $\tau_K$ realized in some $M \in K$, 
a quantifier free $n$-type in $\tau_\Phi$
\sn
\item[$(b)'$]   in clause (b) of part (1), the quantifier free type
which $\langle a_{t_0},\dotsc,a_{t_{n-1}}\rangle$ realizes in $M$ is
$\Phi(\tp_{\qf}(\langle t_0,\dotsc,t_{n-1}\rangle,\emptyset,I))$ for $n <
\omega,t_0,\dotsc,t_{n-1} \in I$.
\end{enumerate}
\end{definition}

\begin{fact}
\label{z8}
1) Let ${\gk}$ be an a.e.c. and $M \in K_{\gk}$ be of cardinality $\ge
\lambda = \beth_{1,1}(\LST_{\gk})$ recalling we may assume
$|\tau_{\gk}| \le \LST_{\gk}$ as usual.

\Then \, there is a $\Phi$ such that $\Phi$ is proper for linear orders and:
\mn
\begin{enumerate}
\item[$(a)$]   $(\alpha) \quad \tau_{\gk} \subseteq \tau_\Phi$,
\sn
\item[${{}}$]   $(\beta) \quad |\tau_\Phi| = \LST_{\gk} + |\tau_{\gk}|$ 
\sn
\item[$(b)$]   for any linear order $I$ the model 
$\EM(I,\Phi)$ has cardinality $|\tau(\Phi)|+ |I|$ and 
we have $\EM_{\tau({\gk})}(I,\Phi) \in K_{\gk}$
\sn
\item[$(c)$]   for any linear orders $I \subseteq J$ we have 
$\EM_{\tau({\gk})}(I,\Phi) \le_{\gk} \EM_{\tau({\gk})}(J,\Phi)$;
  moreover, if $M \subseteq \EM(J,\Phi)$ then $(M \rest \tau_{\gk})
  \le_{\gk} \EM_{\tau(\gk)}(J,\Phi)$
\sn
\item[$(d)$]   for every finite linear order $I$, the model
$\EM_{\tau({\gk})}(I,\Phi)$ can be $\le_{\gk}$-embedded into $M$.
\end{enumerate}
\mn
1A) Moreover, assume in (1) also $\lambda = \beth_{1,1}(\kappa),\kappa
\ge \LST_{\gk} + |\tau_{\gk}|$ so not necessarily assuming $\LST_{\gk}
\ge |\tau_{\gk}|,M^+$ is an expansion of $M$ with
$\tau(M^+)$ of cardinality $\le \kappa$ and $b_\alpha \in M$ for
$\alpha < \lambda$ are pairwise distinct.  \Then \, there is $\Phi$
proper for linear orders such that:
\mn
\begin{enumerate}
\item[$(a)$]  $(\alpha) \quad \tau(M^+) \subseteq \tau_\Phi$ hence
  $\tau(\gk) \subseteq \tau_\Phi$
\sn
\item[${{}}$]  $(\beta) \quad \tau_\Phi$ has cardinality $\kappa$
\sn
\item[$(b),(c)$]   has in part (1)
\sn
\item[$(d)$]  if $I$ is a finite linear order and $t_0 <_I \ldots <_I
  t_{n-1}$ list its elements and $M_I = \EM(I,\Phi)$ with skeleton
  $\langle a_{t_i}:t \in I\rangle$, then for some ordinals $\alpha_0 <
  \ldots < \alpha_{n-1} < \lambda$ there is an embedding of $M_I$ into
  $M^+$ mapping $a_{t_\ell}$ to $b_{\alpha_\ell}$ for $\ell < n$.
\end{enumerate}
\mn
2) If $\LST_{\gk} < |\tau_{\gk}|$ and there is $M \in K_{\gk}$ 
of cardinality $\ge \beth_{1,1}(2^{\LST_{\gk}})$, \then \,
 there is $\Phi \in \Upsilon^{\oor}_{\LST({\gk})+|\tau(\Phi)|} [{\gk}]$ 
such that $\EM(I,\Phi)$ has cardinality $\le \LST_{\gk}$ for $I$
finite and $\tau_\Phi \backslash \tau(M)$ has cardinality
$\LST_{\gk}$.  
Note that ${\cE}$ has $\le 2^{\LST_{\gk}}$ equivalence classes 
where ${\cE} = \{(P_1,P_2):P_1,P_2 \in \tau_\Phi$ and
$P^{\EM(I,\Phi)}_1 = P^{\EM(I,\Phi)}_2$ for every linear order $I\}$
hence above ``$\ge \beth_{1,1}(2^{\LST(\gk)})$" suffice.

\noindent
3) We can combine parts (1A) and (2).  Also in both cases 
having a model of cardinality $\ge \beth_\alpha$ for every
$\alpha < (2^{\LST({\gk})+|\tau({\gk})|})^+$ suffice in parts (1),(1A) and
for every $\alpha < \beth_2(\LST_{\gk})^+$ suffice in part (2). 
\end{fact}

\noindent
We add
\begin{claim}
\label{z9}
For every cardinal $\mu$ and strong limit $\chi \le \mu$ there is a
dense $\kappa$-saturated linear order $I=I_\mu$ of cardinality $\mu$ such that:
\mn
\begin{enumerate}
\item[$(*)$]  if $\partial < \mu,2^\theta \le \mu$ \then \,
\sn
\item[$(*)_{I,\chi,\partial,\theta}$]  we have  $2^\theta \le \mu$ and
  $\theta < \partial = \cf(\partial)$ and $(A) \Rightarrow (B)$ where:
\sn
\begin{enumerate}
\item[$(A)$]  $(a) \quad I_0 \subseteq I$
\sn
\item[${{}}$]  $(b) \quad I_0$ has cardinality $\le \theta$
\sn
\item[${{}}$]  $(c) \quad I_1$ is a linear order extending $I_0$
\sn
\item[${{}}$]  $(d) \quad u_n \subseteq u_{n+1} \subseteq
  \theta = \bigcup\limits_{n} u_n$
\sn
\item[${{}}$]  $(e) \quad \bar t^1_\alpha \in {}^\theta(I_1)$
  for $\alpha < \partial$ and $\langle \bar t^1_\alpha:\alpha <
  \partial \rangle$ is an indiscernible sequence

\hskip25pt  in $I_1$ over $I_0$ (for quantifier free formulas)
\sn
\item[${{}}$]  $(f) \quad$ for every $n,I_{1,n} = I_1 \rest
  (\{t^1_{\alpha,i}i \in u_n,\alpha < \partial\} \cup I_0)$ is
  embeddable

\hskip25pt  into $I$ over $I_0$
\sn
\item[$(B)$]  there is $\langle \bar t_\alpha:\alpha < \mu\rangle$
  such that
\sn
\item[${{}}$]  $(a) \quad \bar t_\alpha \in {}^\theta I$
\sn
\item[${{}}$]  $(b) \quad \langle \bar t_\alpha:\alpha < \mu\rangle$
  is an indiscernible sequence over $I_0$ into $I$ 

\hskip5pt (for quantifier free formulas)
\sn
\item[${{}}$]  $(c) \quad$ the quantifier free type of $\bar t_0 \char
  94 \ldots \char 94 \bar t_n$ over $I_0$ in $I$ is equal to the

\hskip25pt  quantifier free type of $\bar t^1_0 \char 94 \ldots \char 94 \bar
  t^1_n$ over $I_0$ in $I_1$ for every $n$
\sn
\item[$(B)^+$]  moreover we can replace $\langle \bar t_\alpha:\alpha
  < \mu\rangle$ by $\langle \bar t_s:s \in I\rangle$.
\end{enumerate}
\end{enumerate}
\end{claim}

\begin{remark}
1) We may consider replacing (A)(e) by 
\mn
\begin{enumerate}
\item[$(e)'$]   $\alpha = \beth_2(\theta)^+,u_n 
\subseteq u_{n+1} \subseteq \theta =
\bigcup\limits_{n} u_n$ and $I_{1,n} =
\{t^1_{\alpha,\varepsilon}:\alpha < \partial,\varepsilon \in u_n\}$
and there is $\bar f = \langle f_\eta:\eta \in \Lambda\rangle$ such
that $f_\eta$ embeds $I_{1,\ell g(\eta)}$ into $I_1$ over $I_0$ and $\nu
\triangleleft \eta \Rightarrow f_\nu \subseteq f_\eta$ where $\Lambda
= \{\eta:\eta$ is a decreasing sequence of ordinals $< \alpha\}$.
\end{enumerate}
\mn
2) Clauses (A)(d),(e) can be weakened to:
\mn
\begin{enumerate}
\item[$\oplus$]  if $i,j < \theta$ \then \, $I_1 \rest
  (\{t^1_{\alpha,i}\alpha =0,1$ and $i < \theta\} \cup I_0)$ can be
  embedded into $I$ over $I_0$.
\end{enumerate}
\mn
But the present form fits our application.
\end{remark}

\begin{PROOF}{\ref{z9}}
First we give a sufficient condition for $(*)_{I,\partial,\theta}$
\mn
\begin{enumerate}
\item[$\boxplus$]  the linear order $I$ satisfies
$(*)_{I,\partial,\theta}$ when: $\chi > \partial = \cf(\partial) >
  \theta$ and
\sn
\begin{enumerate}
\item[$(a)$]  $I$ is a linear order of cardinality $\mu$
\sn
\item[$(b)$]  if $I_0 \subseteq I,|I_0| \le \theta$ \then \, the set
  $I^+_0 = \{t \in I:t \notin I_0$ and there is no $t' \in I
  \backslash I_0 \backslash \{t\}$ realizing the same cut of $I_0$ in
  $I\}$ has cardinality $< \partial$, so if $\partial = (2^\theta)^+$
  this holds
\sn
\item[$(c)$]  if $a <_I b$ then $I$ is embeddable into $(a,b)_I$
\sn
\item[$(d)$]  every linear order of cardinality $\le \theta$ is
  embeddable into $I$
\sn
\item[$(e)$]  in $I$ there is a decreasing sequence of length $\mu$
  and an increasing sequence of length $\mu$
\sn
\item[$(f)$]  to get $(B)^+$ we need:  if $J$ is a linear order of
  cardinality $\le \theta$ then we can embed $I \times J$ (ordered
  lexicographically into $I$).
\end{enumerate}
\end{enumerate}
\mn
It is obvious that there is such linear order.  It is also each that
if $I$ satisfies (a)-(d) then $(*)_{I,\partial,\theta}$.
\end{PROOF}
\newpage

\section {More on Templates} \label{More}

Why do we need $\Upsilon^{\sor}_\kappa[M,\gk]$?  Remember that such
$\Phi$'s are
witnesses to $M$ having $\le_{\gk}$-extensions in every $\mu >
\LST_{\gk} + \|M\|$ so proving existence is a major theme here.
First, why do we need below $\Upsilon^{\sor}_\kappa$?  Because 
``$\Upsilon^{\sor}_\kappa[M,\gk] \ne \emptyset"$ is equivalent to $M$
as not $\le_{\gk}$-maximal; moreover has $\le_{\gk}$-extensions of
arbitrarily large cardinality so proving this for every $M \in
K^{\gk}_\lambda$ indicates ``$\gk$ is nice, at least in $\lambda$".
Second, why do we need various partial orders on
$\Upsilon^{\sor}_\kappa[M,\gk]$'s?

In a major proof here to build $\Phi \in \Upsilon^{\sor}_\kappa[M,\gk]$ we
use $\le_{\gk}$-increasing $M_n$ with union $M$ and try to choose $\Phi_n
\in \Upsilon^{\sor}_\kappa[M_n,\gk]$ increasing with $n$.  For this we
assume $\|M_n\| = \lambda_n,\lambda_n \ll \lambda_{n+1}$ and we use an
induction hypothesis that $\Phi_n$ has a say $\lambda_{n+5}$-witness in $M$.

Of course, it is nice if $\EM_{\tau(\gk)}(\lambda_{n+5},\Phi_n)$ is
$\le_{\gk}$-embeddable into $M$ over $M_n$ but for this we do not have
strong enough existence theorem.  To fine 
tune this and having a limit $(\Phi \in
\Upsilon^{\sor}_\kappa[M,\gk])$ we need some orders.
\bigskip

\begin{definition}
\label{z12}
For $\gk$ an a.e.c. and $M \in K_{\gk}$ let $\gk_M = \gk[M]$ be the
   following a.e.c.:
\mn
\begin{enumerate}
\item[$(a)$]  vocabulary $\tau_{\gk} \cup \{c_a:a \in M\}$ where the $c_a$'s
are pairwise distinct new individual constants
\sn 
\item[$(b)$]  $N \in K_{\gk_M}$ \Iff \, $N \rest \tau_{\gk} \in
  K_{\gk}$ and $a \mapsto c^N_a$ is a $\le_{\gk}$-embedding of $M$
  into $N \rest \tau_{\gk}$;
\sn
\item[$(c)$]  $N_1 \le_{\gk_M} N_2$ \Iff \,
\sn
\begin{enumerate}
\item[$(\alpha)$]  $N_1,N_2$ are $\tau_{\gk_M}$-models from $K_{\gk_\mu}$
\sn
\item[$(\beta)$]  $N_1 \subseteq N_2$
\sn
\item[$(\gamma)$]  $(N_1 \rest \tau_{\gk}) \le_{\gk} (N_2 \rest
  \tau_{\gk})$.
\end{enumerate}
\end{enumerate}
\end{definition}

\begin{definition}
\label{z14}
1) We call $N \in K_{\gk_M}$ standard \when \, $M \le_{\gk} N \rest
   \tau_{\gk}$ and $a \in M \Rightarrow c^N_a = a$.

\noindent
2) If $N^1 \in K_{\gk_M}$ is standard and $N^0 = N^1 \rest \tau_{\gk}$
   \then \, we write $N^1 = N^0_{[M]}$.

\noindent
3) We call $\Phi \in \Upsilon^{\oor}_{\gk}$ standard \when \, $M =
\EM_{\tau(\gk)}(\emptyset,\Phi)$ implies $N \le_{\gk} M \rest \tau_{\gk}$ when
   $N$ is the submodel\footnote{Note that we have not said ``$\Phi \in
   \Upsilon^{\oor}_{\gk[N]}"$ but by renaming this follows.}
 of $M \rest \tau_{\gk}$ with universe
$\{c^M:c \in \tau(\Phi)$ an individual constant$\}$.  We call $\Phi$ fully
   standard when above $N = M \rest \tau_{\gk}$.

\noindent
4) Let $\Upsilon^{\sor}_\kappa[\gk]$ be the class of standard $\Phi
\in \Upsilon^{\oor}_\kappa[\gk]$.

\noindent
5) For $M \in K_{\gk}$ let 
$\Upsilon^{\sor}_\kappa[\gk_M]$ be the class of $\kappa$-standard
$\Phi \in \Upsilon^{\oor}_\kappa[\gk_M]$ which\footnote{So though such
  $\Phi$ belongs to $\Upsilon^{\oor}_\kappa[\gk]$, being standard for
  $\Upsilon^{\sor}_\kappa[\gk_M]$ is a different demand than being
  standard for $\Upsilon^{\oor}_\kappa[\gk]$ as for the latter
possibly $\{c_a:a \in M\} \subsetneqq \{c \in \tau_\Phi:c$ an individual
constant$\}$.}
means:
\mn
\begin{enumerate}
\item[$(a)$]  letting $\kappa_1 = \kappa + \|M\|$, we have
$\Phi \in \Upsilon^{\sor}_{\kappa_1}[\gk]$
\sn
\item[$(b)$]  $\{c_a:a \in M\} = \{c \in \tau(\Phi):c$ an 
individual constant$\}$. 
\sn
\item[$(c)$]  $N = \EM(\emptyset,\Phi) \Rightarrow |N| =
\{c^N:c \in \tau_\Phi\}$
\sn
\item[$(d)$]  $\tau'_\Phi := \tau_{\gk} \backslash \{c \in \tau_\Phi$ is an
individual constant$\}$ has cardinality $\le \kappa$
\sn
\item[$(e)$]  if $N = \EM(I,\Phi)$ and $N_1$ is a submodel of $N \rest
\tau'_\Phi$ then $N_1 \rest \tau_{\gk} \le_{\gk} N \rest \tau_{\gk}$.
\end{enumerate}
\mn
5A) We may omit $\kappa$ in part (5) when $\kappa = \LST_{\gk} + 
|\tau_{\gk}|$.  We may write $\Upsilon^{\sor}_\kappa[M,\gk]$ instead
of $\Upsilon^{\sor}_\kappa[\gk_M]$ for some $\kappa$, useful 
when $\gk$ is not clear from the context.
\end{definition}

\begin{observation}
\label{z15}
1) If $\Phi \in \Upsilon^{\sor}_\kappa[\gk,M]$ \then \, $\Phi \in
   \Upsilon^{\oor}_{\kappa +\|M\|}[\gk]$ but not necessarily the
   inverse.

\noindent
2) If $\Phi \in \Upsilon^{\sor}_\kappa[\gk,M]$ \then \, $\Phi$ is a
fully standard member of $\Upsilon^{\oor}_\kappa[\gk_M]$.
\end{observation}

\begin{claim}
\label{z17}
Assume $\gk$ is an a.e.c. and $M \in K_{\gk}$ and $\gk_1 = \gk_M$
\then \,:
\mn
\begin{enumerate}
\item[$(a)$]  $\gk_1$ is an a.e.c.
\sn 
\item[$(b)$]  $\LST_{\gk_1} = \LST_{\gk} + \|M\|$
\sn
\item[$(c)$]  applying \ref{z8} to $\gk_1$, we can add $``\Phi \in
\Upsilon^{\sor}_\kappa[\gk_M]"$.
\end{enumerate}
\end{claim}

\begin{PROOF}{\ref{z17}}
Straightforward.
\end{PROOF}

\begin{definition}
\label{z18}
Assume $J$ is a linear order of cardinality $\lambda$ and $\lambda
\rightarrow (\mu)^n_\theta$.  We define the ideal $\cI =
\ER^n_{J,\mu,\theta}$ on the set $[J]^\mu$ by:
\mn
\begin{enumerate}
\item[$\bullet$]  $\cS \subseteq [J]^\mu$ belongs to $\cI$ \Iff \, for
  some $\bold c:[J]^{\le n} \rightarrow \theta$ there is no $s \in
  \cS$ such that $\bold c \rest [s]^n$ is constant.
\end{enumerate}
\end{definition}

\begin{observation}
\label{z8d}
1) If $|J| = \lambda$ and $\lambda \rightarrow (\mu)^n_\theta$ \then \,
   $\ER^n_{J,\mu,\theta}$ is indeed an ideal, i.e. $J \notin
   \ER^n_{J,\mu,\theta}$.

\noindent
2) If $\theta = \theta^{< \kappa}$ \then \, this ideal is $\kappa$-complete.
\end{observation}

\begin{definition}
\label{z19}
1) For vocabularies $\tau_1,\tau_2$ we say that $\bold h$ is an
   isomorphism from $\tau_1$ onto $\tau_2$ \when \, $\bold h$ is a
   one-to-one function from the non-logical symbols of $\tau_1$ (= the
   predicates and function symbols) onto those of $\tau_2$ such that:
\mn
\begin{enumerate}
\item[$(a)$]  if $P \in \tau_1$ is a predicate \then \, $\bold h(P)$
  is a predicate of $\tau_2$ and $\arity_{\tau_1}(P) =
  \arity_{\tau_2}(\bold h(P))$
\sn
\item[$(b)$]  if $F \in \tau_1$ is a function symbol\footnote{this
 includes individual constants}
 \then \, $\bold h(F)$
is a function symbol of $\tau_2$ and $\arity_{\tau_1}(F) =
\arity_{\tau_2}(\bold h(F))$.
\end{enumerate}
\mn
2) If $\bold h$ is an isomorphism from the vocabulary $\tau_1$ onto
the vocabulary $\tau_1$ and $M_1$ is a $\tau_1$-model \then \,
$M^{[\bold h]}_1$ is the unique $M_2$ such that:
\mn
\begin{enumerate}
\item[$(a)$]  $M_2$ is a $\tau_2$-model
\sn
\item[$(b)$]  $|M_2| = |M_1|$
\sn
\item[$(c)$]  $P^{M_2}_2 = P^{M_1}_1$ when $P_1 \in \tau_1$ is a
  predicate and $P_2 = \bold h(P_1)$
\sn
\item[$(d)$]  $F^{M_2}_2 = F^{M_1}_1$ when $F_1 \in \tau_1$ is a
  function symbol and $F_2 = \bold h(F_1)$.
\end{enumerate}
\mn
3) We say $\bold h$ is an isomorphism from $\tau_1$ onto $\tau_2$ over
$\tau$ \when \, $\tau \subseteq \tau_1 \cap \tau_2,\bold h$ is an
isomorphism from $\tau_1$ onto $\tau_2$ and $\bold h \rest \tau$ is
the identity.

\noindent
4) If $\Phi_1 \in \Upsilon^{\oor}_\kappa$ and $\bold h$ is an
   isomorphism from the vocabulary $\tau_1 := \tau(\Phi)$ onto the
   vocabulary $\tau_2$ \then \, $\Phi^{[\bold h]}$ is the unique $\Phi_2
   \in \Upsilon^{\oor}_\kappa$ such that: if $I$ is a linear order,
   $M_1 = \EM(I,\Phi_1)$ with skeleton $\langle a_t:t \in I\rangle$
   then $M^{[\bold h]}_1$ is the model $(\EM(I,\Phi_2))^{[\bold h]}$
   with the same skeleton.
\end{definition}

\begin{observation}
\label{z20}
1) In \ref{z19}(2), $M_2 = M^{[\bold h]}_1$ is indeed a
   $\tau_2$-model.  If in addition $\bold h$ is over 
$\tau$ (i.e. $\tau \subseteq \tau_1 \cap \tau_2$ and 
$\bold h \rest \tau = \id_\tau$) \then \, $M_1 \rest \tau = M_2 \rest
\tau$.

\noindent
2) In \ref{z19}(4), indeed $\Phi_2 \in \Upsilon^{\oor}_\kappa$.

\noindent
3) If $\bold h$ is an isomorphism from $\tau_1$ onto $\tau_2$ over
$\tau_{\gk}$ so $\tau_{\gk} \subseteq \tau_1 \cap \tau_2$ and $\Phi_1
   \in \Upsilon^{\oor}_\kappa[\gk],\tau_1 = \tau(\Phi_1)$ 
\then \, $\Phi_2 = \Phi^{[\bold h]}_1$ belongs to 
$\Upsilon^{\oor}_\kappa[\gk]$.

\noindent
4) In part (3) if in addition $M \in K_{\gk}$ and $\Phi_1 \in
\Upsilon^{\sor}_\kappa[M,\gk]$ and $a \in M \Rightarrow \bold h(c_a) 
= c_a$ \then \, $\Phi_2 = \Phi^{[\bold h]}_1$ belongs to
   $\Upsilon^{\sor}_\kappa[M,\gk]$. 
\end{observation}

\begin{PROOF}{\ref{z20}}
Straightforward.
\end{PROOF}

\noindent
Next we recall the partial orders 
$\le^1_\kappa,\le^2_\kappa$ and define an equivalence
relation and some quasi-orders on $\Upsilon^{\oor}_\kappa[\gk]$.
\begin{definition}
\label{z21}
Fixing $\gk$, we define partial orders
$\le^\oplus_\kappa = \le^1_\kappa = \le^1_{\gk,\kappa}$ 
and $\le^\otimes_\kappa = \le^2_\kappa = \le^2_{\gk,\kappa}$ on 
$\varUpsilon^{\oor}_\kappa[\gk]$ (for $\kappa \ge \LST_{\gk}$): 

\noindent
1)  $\Psi_1 \le^\oplus_\kappa \Psi_2$ \If \, $\tau(\Psi_1) 
\subseteq \tau (\Psi_2)$ and $\EM_{\tau({\gk})}(I,\Psi_1) 
\le_{\gk} \EM_{\tau({\gk})}(I,\Psi_2)$ and $\EM(I,\Psi_1) = 
\EM_{\tau(\Psi_1)}(I,\Psi_1) \subseteq \EM_{\tau(\Psi_1)}(I,\Psi_2)$
for any linear order $I$ (so, of course, same $a_t$'s, etc.). 

\noindent
Again for $\kappa = \LST_{\gk} + |\tau_{\gk}|$ we may drop the $\kappa$. 

\noindent
2) For $\Phi_1,\Phi_2 \in \Upsilon^{\oor}_\kappa$, we say 
$\Phi_2$ is an inessential extension of $\Phi_1$ and write 
$\Phi_1 \le^{\ie}_\kappa \Phi_2$ \If \, $\Phi_1 \le^\oplus_\kappa
\Phi_2$ and for every linear order $I$, we have

\[
\EM_{\tau({\gk})}(I,\Phi_1) = \EM_{\tau({\gk})}(I,\Phi_2).
\]

\[
\text{(note: there may be more function symbols in } \tau(\Phi_2)!)
\]

\mn
2A) We define the two-place relation $\bold E^{\ae}$ on 
$\Upsilon^{\oor}_{\gk}$ as follows $\Phi_1 \bold E^{\ae} \Phi_2$ \Iff
\, $\tau(\Phi_1) = \tau(\Phi_2)$ and for some unary function symbol $F
\in \tau(\Phi_1)$ or $F$ is just a (finite)
composition\footnote{but abusing our notation we may still write $F \in
\tau_\Phi$} of such function symbols, if $M = 
\EM(I,\Phi_1)$ with skeleton $\langle a^1_t:t \in
I\rangle$ and we let $a^2_t = F^M(a^1_2)$ for $t \in I$ then:
\mn
\begin{enumerate}
\item[$\bullet$]  $F^M(a^2_t) = a^1_t$
\sn
\item[$\bullet$]  $M$ is $\EM(I,\Phi_2)$ with skeleton $\langle a^2_t:t \in
I\rangle$;
\end{enumerate}
\mn
``$\ae$" stands for almost equal.

\noindent
2B) Above we say $\Phi_2 \bold E^{\ae} \Phi_2$ is witnessed by $F$.

\noindent
2C) We define the two-place relation $\bold E^{\ie}_\kappa$ on
$\Upsilon^{\oor}_{\gk}$ by: $\Phi_1 \bold E^{\ie}_\kappa \Phi_2$ iff
for some $\Phi_3,\Phi_1 \le^{\ie}_\kappa \Phi_3$ and $\Phi_2
\le^{\ie}_K \Phi_3$.

\noindent
2D) We define a two-place relation $\bold E^{\ai}_\kappa$ on
$\Upsilon^{\oor}_\kappa[\gk]$ by $\Phi_1 \bold E^{\ai}_\kappa \Phi_3$ iff
for some $\Phi_2 \in \Upsilon^{\ai}_\kappa[\gk]$ we have
$\Phi_1 \bold E^{\ae}_\kappa \Phi_2$ and $\Phi_2 \bold E^{\ie}_K
\Phi_3$.

\noindent
3) Let $\Upsilon^{\lin}_\kappa$ be the class of 
$\Psi$ proper for linear order and producing linear orders, that is, such that:
\begin{enumerate}
\item[$(a)$]    $\tau(\Psi)$ has cardinality $\le \kappa$,
\sn
\item[$(b)$]   $\EM_{\{<\}}(I,\Psi)$ is a linear order 
which is an extension of $I$ which means $s <_I t \Rightarrow
\EM(I,\Psi) \models ``a_s < a_t"$; in fact we can have

$[t \in  I \Rightarrow a_t = t]$.
\end{enumerate}
\mn
4) $\Phi_1 \le^\otimes_\kappa \Phi_2$ \Iff \, there is $\Psi$ such that:
\begin{enumerate}
\item[$(a)$]   $\Psi \in \Upsilon^{\lin}_\kappa$
\sn
\item[$(b)$]   $\Phi_\ell \in \Upsilon^{\oor}_\kappa$ for $\ell=1,2$
\sn
\item[$(c)$]   $\Phi'_2 \le^{\ie}_\kappa \Phi_2$ where
$\Phi'_2 = \Psi \circ \Phi_1$, i.e.
\end{enumerate}

\[
\EM_{\tau(\Phi_1)}(I,\Phi'_2) = \EM(\EM_{\{<\}}(I,\Psi),\Phi_1).
\]

\mn
(So we allow further expansion by functions definable from earlier ones
(composition or even definition by cases), as long as the number is 
$\le \kappa$). 
\end{definition}

\noindent
It is not a real loss to restrict ourselves to standard $\Phi$ because
\begin{claim}
\label{z22}
1) For every $\Phi_1 \in \Upsilon^{\oor}_\kappa[\gk]$ there is a standard
$\Phi_2 \in \Upsilon^{\oor}_\kappa[\gk]$ such that $\Phi_1
\le^{\ie}_\kappa \Phi_2$; moreover $M = \EM(\emptyset,\Phi_2)
\Rightarrow |M| = \{c^M:c \in \tau(\Phi_2)$ an individual
constant$\}$, that is $\Phi_2$ is fully standard.

\noindent
2) Assume $\Phi_1 \in \Upsilon^{\oor}_\kappa[\gk],F \in \tau(\Phi)$ is a
   unary function symbol such that $M = \EM(I,\Phi_1) \wedge t \in I
   \Rightarrow F^M(F^M(a_t))=a_t$. \Then \, for a unique
   $\Phi_2,\Phi_1 \bold E^{\ae} \Phi_2$ as witnessed by $F$ and
   $\Phi_1 \in \Upsilon^{\sor}_\kappa[\gk_M] \Leftrightarrow \Phi_2
   \in \Upsilon^{\sor}_\kappa[\gk_M]$.

\noindent
3) $\bold E^x_\kappa$ is an equivalence relation on
   $\Upsilon^{\oor}_K[\gk]$ for $\lambda \in \{\ae,\ie,\ai\}$ all
   refining $\bold E^{\ai}_\kappa$.
\end{claim}

\begin{PROOF}{\ref{z22}}
Obvious.
\end{PROOF}

\begin{observation}
\label{z23}
Let $\ell=1,2$.

\noindent
1) The relation $\le^\ell_\kappa$ is a partial order on
   $\Upsilon^{\oor}_\kappa[\gk]$.

\noindent
2) If $\langle \Phi_\alpha:\alpha < \delta\rangle$ is
   $\le^\ell_\kappa$-increasing with $\delta$ a limit ordinal $<
   \kappa^+$ \then \, $\bigcup\limits_{\alpha < \delta}
   \Phi_\alpha$ naturally defined is a $\le^\ell_\kappa$-lub.

\noindent
3) $\bold E^{\ae}$ is an equivalence relation on $\Upsilon^{\oor}$.

\noindent
4) If $\Upsilon^{\oor}_{\kappa_1}[\gk] \subseteq
   \Upsilon^{\oor}_{\kappa_2}[\gk]$ \then \, $\kappa_1 \le
\kappa_2$.  If $\kappa_1 \le \kappa_2$ and 
$\iota \in \{1,2\}$ and $\Phi,\Psi \in
\Upsilon^{\oor}_{\kappa_1}$ \then \, $[\Phi \le^\iota_{\kappa_1} \Psi
   \Leftrightarrow \Phi \le^\iota_{\kappa_2} \Psi)$.

\noindent
5) Similarly for $\Upsilon^{\sor}_\kappa[\gk_M]$ defined in \ref{z14}(5).
\end{observation}

\begin{definition}
\label{z24}
1) For $\kappa \ge \LST_{\gk} + |\tau_{\gk}|$, we define
$\le^\odot_\kappa = \le^3_\kappa$, in full $\le^3_{\gk,\kappa}$, 
a two-place relation on $\Upsilon^{\sor}_\kappa[\gk]$, recalling
Definition \ref{z14}(5) as follows:

Let $\Phi_1 \le^3_\kappa \Phi_2$ mean that: for every linear order
$I_1$ there are a linear order $I_2$ and $\le_{\gk}$-embedding $h$ of
$\EM_{\tau(\gk)}(I_1,\Phi_1)$ into $\EM_{\tau(\gk)}(I_2,\Phi_2)$,
moreover every individual constant $c$ of $\tau(\Phi_1)$ is 
an individual constant of $\tau(\Phi_2)$ and 
$h(c^{\EM(I_1,\Phi_1)}) = c^{\EM(I_2,\Phi_2)}$.

\noindent
2) We define $\le^4_\kappa = \le^4_{\gk,\kappa}$; a two-place
relation on $\Upsilon^{\sor}_\kappa[\gk]$ as follows.

Let $\Phi_1 \le^4_\kappa \Phi_2$ mean that: for some $F$ we have:
\mn
\begin{enumerate}
\item[$(a)$]  $\Phi_1,\Phi_2 \in \Upsilon^{\sor}_\kappa[\gk]$
\sn
\item[$(b)$]  $\bullet \quad \tau(\Phi_1) \subseteq \tau(\Phi_2)$
\sn
\item[${{}}$]  $\bullet \quad F \in \tau(\Phi_2)$ is a unary function
symbol or as in \ref{z21}(2A) 
\sn
\item[$(c)$]  if $I$ is a linear order and 
$M_2 = \EM(I,\Phi_2)$ with skeleton $\langle a^2_s:s
  \in I\rangle$ \then \, there is $M_1 = \EM(I,\Phi_1)$ with skeleton
  $\langle a^1_s:s \in I\rangle$ such that
\sn
\begin{enumerate}
\item[$\bullet$]  $a^1_s = F^{M_2}(a^2_s)$ for $s \in I$
\sn
\item[$\bullet$]  $a^2_s = F^{M_2}(a^1_s)$ for $s \in I$
\sn
\item[$\bullet$]  $M_1 \subseteq M_2 \rest \tau_{\Phi_1}$ so
  $\tau(\Phi_1) \subseteq \tau(\Phi_2)$
\sn
\item[$\bullet$]  $(M_1 \rest \tau_{\gk}) \le_{\gk} (M_2 \rest
\tau_{\gk})$
\sn
\item[$\bullet$]  $c^{M_1} = c^{M_2}$ \when \, $c \in \tau(\Phi_1)$ is
an individual constant.
\end{enumerate}
\end{enumerate}
\end{definition}

\begin{remark}
\label{z25}
So $\le^4_\kappa$ is like $\le^1_\kappa$ but we demand 
less as $a^1_s = a^2_s$ is weakened by using
the function symbol $F$.
\end{remark}

\begin{claim}
\label{z26}
1) $\le^3_\kappa$ is a partial order on $\Upsilon^{\sor}_\kappa[\gk]$
 as well as $\le^4_\kappa$; also for $\Phi_1,\Phi_2 \in
 \Upsilon^{\sor}_\kappa[\gk]$ and $\ell=1,2,4$ we have $\Phi_1 \le^2
 \Phi_2 \Rightarrow \Phi_1 \le^1 \Phi_2 \Rightarrow \Phi_1 \le^4
 \Phi_2 \Rightarrow \Phi_1 \le^3 \Phi_3$.

\noindent
2) Assume $\Phi_1,\Phi_2 \in \Upsilon^{\sor}_\kappa[\gk]$ have the
   same individual constants.  Then
$\Phi_1 \le^3_\kappa \Phi_2$ \Iff \, as in \ref{z24}(1) restricting
ourselves to $I = \beth_{1,1}(\kappa)$ \Iff \,
$\Phi_1,\Phi_2 \in \Upsilon^{\sor}_\kappa[\gk]$ and for some $F$ and $\Phi'_1,
\Phi'_2 \in \Upsilon^{\sor}_\kappa[\gk]$ and we have $\Phi_1 \le^4_\kappa 
\Phi'_1$ witnessed by $F$ and $\Phi'_1 \bold E^{\ae} \Phi'_2$
witnessed by $F$ and for some $\tau_*,\bold h$ we have 
$\tau(\gk) \subseteq \tau_* \subseteq \tau(\Phi'_1),h$ 
is an isomorphism from $\tau(\Phi_2)$ onto $\tau_*$ over $\tau(\gk)
\cup \{c:c \in \tau(\Phi_1)\}$ and $\Phi^{[\bold h]}_2 
\le^{\ie}_\kappa \Phi'_2$ \Iff \, for some $\Phi' \in
   \Upsilon^{\sor}_\kappa[\gk]$ we have $\Phi_1 \le^3 \Phi'$ and
   $\Phi' \bold E^{\ai}_\kappa \Phi$, see \ref{z21}(2).

\noindent
3) If $\Phi_n \in \Upsilon^{\sor}_\kappa[\gk]$ and $\Phi_n \le^3_\kappa
   \Phi_{n+1}$ \then \, there is $\Phi_\omega \in
   \Upsilon_\kappa[\gk]$ such that $n < \omega \Rightarrow \Phi_n
   \le^3_\kappa \Phi$; moreover, $\EM_{\tau(\gk)}(\emptyset,\Phi)$ is
   the union of the $\le_{\gk}$-increasing sequence $\langle
   \EM_{\tau(\gk)}(\emptyset,\Phi_n):n < \omega\rangle$.

\noindent
4) Similarly for $\le^4_\kappa$.
\end{claim}

\begin{PROOF}{\ref{z26}}
1) Obvious.

\noindent
2) \underline{First clause implies second clause}

Holds trivially.
\smallskip

\noindent
\underline{Second clause implies the third clause}

Let $I_1 = (\lambda,<),\lambda$ large enough, e.g. $\lambda =
\beth_{1,1}(\kappa)$.
Let $M_1 = \EM(I_1,\Phi_1)$ be with skeleton $\langle a^1_t:t \in
I_1\rangle$.  As $\Phi_1 \le^3_\kappa \Phi_2$, there is a linear order
$I_2$ and $M_2 = \EM(I_2,\Phi_2)$ with skeleton $\langle a^2_t:t \in
I_2\rangle$ and $\le_{\gk}$-embedding $f$ from $M_1 \rest \tau(\gk)$
into $M_2 \rest \tau(\gk)$ such that $c \in \tau(\Phi_1) \Rightarrow
c \in \tau(\Phi_2) \wedge f(c^{M_1}) = 
c^{M_2}$; so by renaming \wilog \, $f \rest
\Sk(\emptyset,M_1)$ is the identity and $|I_2| > \lambda$.  
As $\|M_2\| > \|M_1\| \ge
\lambda > \kappa \ge |\tau(M_2)|$, clearly we can find pairwise distinct
$t_\alpha \in I_2$ for $\alpha < \lambda$ such that
$\{a^2_{t_\alpha}:\alpha < \lambda\} \cap \{f(a^1_\alpha):\alpha <
\lambda\} = \emptyset"$.

Let $\tau_1 = \tau(\Phi_1)$ and\footnote{The reason is that there may
  be a symbol in $\tau(\Phi_2) \cap \tau(\Phi_c)$ but not from
  $\tau(\gk_1) \cup \{c:c \in \tau(\Phi_1)\}$. We eliminate this
  ``accidental equality".  Only now $\tau_3 \cup \tau_1$ ``makes sense".}
 let the pair $(\bold h,\tau_3)$ be such that: $\bold h$ 
is an isomorphism from the vocabulary $\tau_2 = 
\tau(\Phi_2)$ onto $\tau_3$ over $\tau(\gk) \cup \{c:c \in \tau(\Phi_1)\}$
such that $\tau_1 \cap \tau_3 = \tau(\gk) \cup \{c:c \in
  \tau(\Phi_2)\}$ and let $M_3 = M^{[\bold h]}_2$, so
$\tau(M_3) = \tau_3,\Phi_3 = \Phi^{[\bold h]}_2$ so $\tau(M_3) = \tau_3 =
\tau(\Phi_3)$ and $M_3$ is an $\EM(I_2,\Phi_3)$ model with skeleton
$\langle a^2_t:t \in I_2 \rangle$.

Let $\tau_4 = \tau_3 \cup \tau_1 \cup \{F,P\}$ with $F$ a one 
place function symbol and $P_\ell,F \notin \tau_3 \cup \tau_1$ and $P_\ell$ 
one place predicates for $\ell=1,2,3,4$.  We define a $\tau_4$-model $M_4$:
\mn
\begin{enumerate}
\item[$\bullet_1$]  it has universe $|M_3|$
\sn
\item[$\bullet_2$]  $F^{M_4}(a^2_{t_\alpha}) =
  f(a^1_\alpha)$ and $F^{M_4}(f(a^1_\alpha)) = a^2_{t_\alpha}$
\sn
\item[$\bullet_3$]  $P^{M_4}_1 =\{a^2_t:t \in I_2\},P^{M_4}_2 =
  \{a^2_t:t \in I_2\},P^{M_4}_3 = \{f(a^1_t):t \in I_1\},P^{M_4}_4 = \Rang(f)$
\sn
\item[$\bullet_4$]  $M_4 \rest \tau_3 = M_3$
\sn
\item[$\bullet_5$]  $f$ embeds $M_1$ into $M_4 \rest \tau_1$.
\end{enumerate}
\mn
Clearly there is no problem to do this and we apply \ref{z8}(1A) with
$M_4 \rest \tau(\gk),M_4,\langle a^2_{t_\alpha}:\alpha <
  \lambda\rangle$, here standing for $M,M^+,\langle b_\alpha:\alpha <
  \lambda\rangle$ there and get $\Phi_4$ standing for $\Phi$ there.
  Now by inspection (see Definition \ref{z24}(2)):
\mn
\begin{enumerate}
\item[$(*)_1$]  $\Phi_1 \le^4_\kappa \Phi_4$
\sn
\item[$(*)_2$]  $\Phi_3 \le^\otimes_\kappa \Phi_4$; moreover $\Phi_3
  \le^{\ie} \Phi_4$.
\end{enumerate}
\mn
We derive $\Phi_5$ from $\Phi_4$ by \ref{z22}(2) using our $F$ so
$\Phi_4 \bold E^{\ae} \Phi_5$.
To show that the third clause of part (2) indeed holds, we just note
that $\Phi'_1,\Phi'_2,\bold h,\tau_*$, there can stand for 
$\Phi_4,\Phi_5,\bold h,\tau_3$ here, so we are done.
\medskip

\noindent
\underline{The third clause implies the first clause}:

So we are given $F$ and 
$\Phi_1,\Phi_2 \in \Upsilon^{\sor}_\kappa[\gk],\Phi'_1,\Phi'_2
\in \Upsilon^{\sor}_\kappa[\gk],\tau_* \subseteq \tau(\Phi'_2)$
including $\tau(\gk)$ and an isomorphism $\bold h$ from $\tau(\Phi_2)$ onto
$\tau_*$ over $\tau_{\gk} \cup \{c:c \in \tau(\Phi_1)\}$
 such that $\Phi_1 \le^4_\kappa \Phi'_2$ witness by $F,\Phi'_1 \bold
E^{\ae} \Phi'_2$ witness by $F$ and $\Phi^{[\bold h]}_2
\le^\otimes_\kappa \Phi'_2$.  

Let $\Psi \in \Upsilon^{\lin}_\kappa$ witness $\Phi^{[\bold h]}_2 
\le^\otimes_\kappa \Phi'_2$; and 
for uniformity of notation we let $\Phi_3 := \Phi'_2$.

We have to prove $\Phi_1 \le^3_\kappa \Phi_2$ so let $I_1$ be a linear
order.

Let $M^*_1 = \EM(I_1,\Phi_1)$ be with skeleton $\langle a^1_t:t \in
I_1\rangle$, let $I_2 = \EM_{\{<\}}(\Psi,I_1)$ so with skeleton
$\langle t:t \in I_1\rangle$.  Let $M_1 \subseteq M_2$ be defined by
$M_\ell = \EM(I_\ell,\Phi_2)$ with skeleton $\bar a^\ell = \{a^2_t:t \in
I_\ell\}$ for $\ell=1,2$ and let $M_3 = \EM(I_1,\Phi'_1)$ be with
skeleton $\{a^3_t:t \in I_1\}$.  

By the choice of $\Psi$ and of $I_2$ \wilog \, $M^{[\bold h]}_2 = M_3 \rest
\tau_*$.

Lastly, there is a unique embedding $f$ of $M^*_1$ into $M_3 \rest
\tau(\Phi_1)$ mapping $a^1_t$ to $F^{M_3}(a^2_t)$ for $t \in I_1$.
Easily $f$ is a $\le_{\gk}$-embedding of $M_1 \rest \tau(\gk)$ into
$M_3 \rest \tau(\gk)$ mapping $c^{M_1}$ to $c^{M_2}$ for $c \in
\tau(\Phi_1)$ and $M_3 \rest \tau(\gk) = M_2 \rest \tau(\gk)$ and $c
\in \tau(\Phi_1) \Rightarrow c \in \tau(\Phi_2) \wedge 
f(c^{M^*_1}) = c^{M_2}$.  

We leave the \underline{fourth} clause to the reader.

\noindent
3) By parts (2) and (4) \underline{or} directly using \ref{z8}(1) and
   the definition of $\le^3_\kappa$.

\noindent
4) So assume that $n < \omega \Rightarrow \Phi_n \le^4_\kappa
   \Phi_{n+1}$ as witnessed by $F_n \in \tau(\Phi_{n+1})$.  For any
   infinite linear order $I$ we can choose $M_n = \EM(I_n,\Phi_n)$
   with skeleton $\langle a^n_t:t \in I\rangle$.  Let $\tau_\omega =
   \cup\{\tau(\Phi_n):n < \omega\}$.  \Wilog \, $M_n \subseteq M_{n+1} 
\rest \tau(\Phi_n),F^{M_{n+1}}_n(a^{n+1}_t) = a^n_t$ and
   $F^{M_{n+1}}_n(a^n_t) = a^{n+1}_t$.  For each $n$ we define $M_{\omega,n} = 
\cup\{M_{n+k} \rest \tau_n:k \in [n,\omega)\}$, so $n_1 <
   n_2 \Rightarrow M_{\omega,n_1} = M_{\omega,n_2} \rest
   \tau(\Phi_{n_1})$.  Hence letting $\tau_\omega =
   \cup\{\tau(\Phi_n):n < \omega\}$ there is a $\tau_\omega$-model
   $M_\omega$ with universe $|M_{\omega,0}|$ such that $M_\omega \rest
   \tau_n = M_{\omega,n}$ for $n < \omega$.  
Now define $\Phi$ by $\Phi(n) = \tp_{\qf}(\langle
   a^0_{t_0},\dotsc,a^0_{t_{n-1}}\rangle,\emptyset,M_\omega)$ whenever
   $t_0 <_I \ldots <_I t_n$.  

Clearly $M_\omega = \EM(I,\Phi)$ with skeleton $\langle a^0_t:t \in
I\rangle$ and $F_{n-1},\dotsc,F_1,F_0$ witness 
$\Phi_n \le^4_\kappa \Phi_\omega$, here we need composition of unary functions.
\end{PROOF}

\begin{claim}
\label{z29}
For $M \in K_{\gk}$ of cardinality $\kappa \ge \LST_{\gk} + |\tau_{\gk}|$ the
following conditions are equivalent:
\mn
\begin{enumerate}
\item[$(a)$]  $\Upsilon^{\oor}_\kappa[\gk_M] \ne \emptyset$
\sn
\item[$(b)$]  for every $\lambda \ge \kappa$ there is $N$ such that $M
  \le_{\gk} N \in K^{\gk}_\lambda$
\sn
\item[$(c)$]  for every $\alpha < (2^\kappa)^+$ there is $N \in
  K^{\gk}_{\ge \beth_\alpha}$ which $\le_{\gk}$-extend $M$
\sn
\item[$(d)$]  there is $\Phi \in \Upsilon^{\oor}_\kappa[\gk_M]$ such that if
$N = \EM(I,\Phi)$ and $N \rest \tau_{\gk_M}$ is standard then $M =
(N \rest \tau_{\gk}) \rest \{c^N:c \in \tau_\Phi$ an individual constant$\}$
\sn
\item[$(e)$]   $\Upsilon^{\sor}_\kappa[\gk_M]$ is non-empty.
\end{enumerate}
\end{claim}

\begin{PROOF}{\ref{z29}}
For (d) note that we can replace an individual constant by a unary
function which is interpreted as being a constant function.
More generally an $n$-place function $F^N$ by functions $F_1,F_2$ \where \,
$\{F_{u,n}:u \subseteq n,h:u \rightarrow M\}$ where
\sn
\begin{enumerate}
\item[$\bullet$]  $F_1$ is a $(n+1)$-place function
\sn
\item[$\bullet$]  if $\bar a = \langle a_\ell:\ell \le n\rangle \in
  {}^{n+1}N \backslash {}^{n+1}M$ then $F_2(\bar a) = F^N(\bar a
  \rest n)$
\sn
\item[$\bullet$]  if $\bar a \in {}^n M$ then $F_1(\bar a) = a_0$
\end{enumerate}
\end{PROOF} 

\begin{claim}
\label{z30}
If (A) then (B) when:
\mn
\begin{enumerate}
\item[$(A)$]  $(a) \quad M_1 \le_{\gk} M_2$
\sn
\item[${{}}$]  $(b) \quad \Phi_1,\Psi_1$ are from
$\Upsilon^{\sor}_\kappa[\gk_{M_1}]$ so are $\kappa$-standard 
\sn
\item[${{}}$]  $(c) \quad \Psi_2 \in
  \Upsilon^{\sor}_\kappa[\gk_{M_2}]$
\sn
\item[${{}}$]  $(d) \quad \Phi_1 \le^4_\kappa \Psi_1$
\sn
\item[${{}}$]  $(e) \quad \Psi_1 \le^1_\kappa \Psi_2$
\sn
\item[${{}}$]  $(f) \quad \{c_a:a \in M_2\} \cap \tau(\Psi_1) =
  \{c_a:a \in M_1\}$
\sn
\item[$(B)$]  there is $\Phi_2$ such that
\sn
\item[${{}}$]  $(a) \quad \Phi_2 \in
  \Upsilon^{\sor}_\kappa[\gk_{M_2}]$ 
\sn
\item[${{}}$]  $(b) \quad \Phi_1 \le^1_\kappa \Phi_2$
\sn
\item[${{}}$]  $(c) \quad \Phi_2 \le^4_\kappa \Psi_2$.
\end{enumerate}
\end{claim}

\begin{PROOF}{\ref{z30}}
Straightforward: let $I$ be an infinite linear order, $M_2 =
\EM(I,\Psi_2)$ be with skeleton $\langle a^2_t:t \in I\rangle$.
Let the unary function symbol $F$ witness $\Phi_1 \le^4_\kappa \Psi_1$
so $F \in \tau(\Psi_1) \subseteq \tau(\Psi_2)$ and let $a^1_t =
F^{M_2}(a^2_t)$.  Clearly $\langle a^1_t:t \in I\rangle$ is
indiscernible for quantifier formulas in $M_2$ and generate it hence
for some $\Phi_2 \in \Upsilon^{\oor}_\kappa$ we have $M_2 =
\EM(I,\Phi_2)$.  Clearly $\Phi_2 \in \Upsilon^{\sor}_\kappa[\gk]$ and
the model $M_2$ is $\EM(I,\Phi_2)$ with skeleton $\langle a^1_t:t \in
I\rangle$.  Also $\Phi_2 \bold E^{\ae} \Phi_2$ hence $\Phi_2
\le^4_\kappa \Psi_2$ and $\Phi_1 \le^{\oplus}_\kappa \Phi_2$ as required.
\end{PROOF}
\bigskip

\centerline {$* \qquad * \qquad *$}
\bigskip

\noindent
The following will be used when applied to a tree of approximations to
embedding of $\EM$-models to a model.  In fact, we use only \ref{z35}
for the case $\cS = \cT \backslash \max(\cT)$, see background in \ref{z37}.
\begin{definition}
\label{z32}
1) We say $\bold i = (\cT,\bar{\bold I}) = (\cT_{\bold i},
\bar{\bold I}_{\bold i})$ is $\pit$ (partially idealized tree) \when \,:
\mn
\begin{enumerate}
\item[$(a)$]  $\cT$ is a tree with $\le \omega$ levels and
\sn
\item[${{}}$]  $\bullet \quad$ for transparency it is a set of finite
  sequences ordred by $\triangleleft$, closed under 

\hskip25pt  initial segments
\sn
\item[${{}}$]  $\bullet \quad$ let $\lev(\eta,\cT) = 
\lev_{\cT}(\eta)$ be the level of $\eta \in \cT$
  in $\cT$, that is 

\hskip25pt  $|\{\nu \in \cT:\nu \triangleleft \eta\}|$
\sn
\item[${{}}$]  $\bullet \quad$ let $\rt_{\cT}$ be the root 
\sn
\item[${{}}$]  $\bullet \quad$ the $n$-level of $\cT$ is the set
  $\{\eta:\lev_{\cT}(n)=n\}$ 

\quad so we have
\sn
\item[${{}}$]  $\bullet \quad \lev_{\cT}(\eta) = 
\ell g(\eta)$ and $\rt_{\cT} = \langle \rangle$
\sn
\item[$(b)$]  $\bold I = \langle \bold I_\eta:\eta \in \cS\rangle$
  where $\cS \subseteq \cT \backslash \max(\cT)$, we may write
$\cS_{\bold i} = \cS$
\sn
\item[$(c)$]  $\bold I_\eta$ is an ideal on $\suc_{\cT}(\eta) :=
  \{\rho:\nu \in \cT,\eta <_{\cT} \rho$ and there is no $\nu \in \cT$
  satisfying $\eta <_{\cT} \nu <_{\cT} \rho\}$ or just an ideal on a set
which $\supseteq \suc_{\cT}(n)$ such that $\suc_{\cT}(\eta) \notin
\bold I_\eta$; we may write $\bold I_{\bold i,\eta}$.
\end{enumerate}
\mn
1A)  If $\bold I_\eta = \{\{s:\eta \char 94 \langle s \rangle \in
X\}:X \in \bold I'_\eta\}$ for some ideal $\bold I'_\eta$ on some set
\then \, abusing 
notation we may write $\bold I'_\eta$ instead of $\bold I_\eta$.

\noindent
2) Let $(\cT_1,\bar{\bold I}_1) \le (\cT_2,\bar{\bold I}_2)$ \when \,
(each is a $\pit$ and):
\mn
\begin{enumerate}
\item[$(a)$]  $\cT_1 \subseteq_{\tr} \cT_2$ which means:
\sn
\begin{enumerate}
\item[$(\alpha)$]   $\eta \in \cT_2 \Rightarrow \eta_1 
\in \cT_1 \wedge \lev(\eta,\cT_2) =
\lev(\eta,\cT_1) \wedge \suc(\eta,\cT_2) \subseteq \suc(\eta,\cT_1)$ 
\sn
\item[$(\beta)$]   $\le_{\cT_1} = <_{\cT_2} \rest \cT_1$
\end{enumerate}
\sn
\item[$(b)$]  $\bar{\bold I}_2 = \bar{\bold I}_1 \rest \cT_2$,
i.e. $\bar{\bold I}_1 \rest \{\eta \in \cT_1:\eta \in \Dom(\bold
I_2)\}$
\sn
\item[$(c)$]  if $\eta \in \cT_2 \backslash \cS_{\bold i_2}$ then
$\suc(\eta,\cT_2) = \suc(\eta,\cT_1)$.
\end{enumerate}
\mn
2A) Let $(\cT_1,\bar{\bold I}_1) \le_{\pr} (\cT_2,\bar{\bold I}_2)$
\when \, (each is a pit and)
\mn
\begin{enumerate}
\item[$(a),(b),(c)$]  as above
\sn
\item[$(d)$]  if $\eta \in \Dom(\bar{\bold I}_1)$ then
  $\suc_{\cT_1}(\eta) \backslash \suc_{\cT_2}(\eta) \in \bold
  I_{1,\eta}$.
\end{enumerate}
\mn
3) We say $(\cT,\bar{\bold I})$ is $\kappa$-complete \when \, every
ideal $\bold I_\eta$ is.

\noindent
4) For $\bold i = (\cT,\bar{\bold I})$ we define $\Dp_{\bold i} =
\Dp_{\cT,\bar{\bold I}}:\cT \rightarrow \Ord 
\cup \{\infty\}$ by (stipulate $\infty + 1 =
\infty$) defining when $\Dp_{\cT,\bar{\bold I}}(\eta) \ge \alpha$
by induction on $\alpha$ as follows:
\mn
\begin{enumerate}
\item[$(a)$]  if $\eta \in \max(\cT)$ \then \,
$D_{\cT,\bold I}(\eta) \ge \alpha$ iff $\alpha = 0$
\sn
\item[$(b)$]  if $\eta \in \cT \backslash \max(\cT)$ and $\eta \in
  \cS_{\bold i} = \Dom(\bar{\bold I})$ \then \, $\Dp_{\cT,\bold I}(\eta)
\ge \alpha$ iff $(\forall \beta < \alpha)(\exists X \subseteq
\suc_{\cT}(\eta))[X \in \bold I^+_\eta \wedge (\forall \nu \in X)
(\Dp_{(\cT,\bar{\bold I)}}(\nu) \ge \beta)]$
\sn
\item[$(c)$]  if $\eta \in \cT \backslash \max(\cT) \backslash
  \cS_{\bold i}$ \then \, $\Dp_{\bold i}(\eta) \ge \alpha$ \Iff \,
  $(\forall \nu)(\nu \in \suc_{\cT}(\eta) \Rightarrow \Dp_{\bold
  i}(\nu) \ge \alpha)$.
\end{enumerate}
\mn
6) If $\bold i = (\cT,\bar{\bold I})$ is a $\pit$ and $\eta \in \cT$
let $\proj(\eta,\bold i) = \proj_{\bold i}(\eta)$ is the sequence $\nu$
of length $\ell g(\eta)$ such that:
\mn
\begin{enumerate}
\item[$\bullet$]  $\ell < \ell g(\eta) \wedge \eta \rest \ell \in
  \Dom(\bar{\bold I}) \Rightarrow \nu(\ell) = -1$
\sn
\item[$\bullet$]  $\ell < \ell g(\eta) \wedge \eta \rest \ell \notin
  \Dom(\bar{\bold I}) \Rightarrow \nu(\ell) = \eta(\ell)$.
\end{enumerate}
\mn
7) For $\bold i = (\cT,\bar{\bold I})$ a $\pit$ let
$\proj(n,\bold i) = \proj_{\bold i}(n) = \{\proj_{\bold i}(\eta):\eta
\in \cT$ has length $n\}$ and $\proj_{\bold i} = \proj(\bold i)$ is
$\{\proj_{\bold i}(\eta):\eta \in \cT\}$.

\noindent
8) If $\bold i_\ell$ is a $\pit$ for $\ell < n$ then
\mn
\begin{enumerate}
\item[$(a)$]  $\prod\limits^{*}_{\ell < \bold n} \cT_{\bold i_\ell}$ is
$\{\bar\eta:\bar\eta = \langle \eta_\ell:\ell < \bold n\rangle$ is
   such that $\ell < \bold n \Rightarrow \eta_\ell \in \cT_{\bold
   i_\ell}$ and moreover for some $n$ called $\lev(\bar\eta)$ we have $(\forall
\ell < \bold n)(\lev_{\cT_{\bold i(\ell)}}(\eta_\ell) = n)\}$. 
\end{enumerate}
\end{definition}

\begin{theorem}
\label{z35}
There are a $\pit \,\bold i_2$ and $\langle c_\eta:\eta \in 
\proj(\bold i_1)\rangle$ such that:
$\bold i_1 \le \bold i_2,\Dp_{\bold i_2}(\rt_{\bold i_2}) \ge
\gamma_2$ and $\eta \in \cT_{\bold i_2} \Rightarrow \bold c(\eta) =
\bold c_{\proj(\eta,\bold i_1)}$ \when \,:
\mn
\begin{enumerate}
\item[$(a)$]   $\bold i_1 = (\cT_1,\bar{\bold I}_1)$ is a $\pit$
\sn
\item[$(b)$]  $\bold i_1$ is $\lambda$-complete $\pit$
\sn
\item[$(c)$]  $2^{\kappa^\theta} < \lambda$ where $\theta =
  |\proj_{\bold i_1}|,\kappa + \theta$ is infinite for
  transparency\footnote{If $\kappa$ and $\theta$ are finite, the
    computations are somewhat different.  Note that $\kappa=0$ is
    impossible and if $\kappa=1$ then $\bold i_2 = \bold i_1$ will do
    so, \wilog \, $\kappa \ge 2$.}
\sn
\item[$(d)$]  $\bold c$ is a colouring of $\cT_1$ by $\le \kappa$ colours
\sn
\item[$(e)$]  $\gamma_1 = \gamma_2 = (2^{(\kappa^\theta)})^+$ or just
\sn
\item[${{}}$]  $(\alpha) \quad \gamma_1 \le \Dp_{\bold i_1}
(\rt_{\bold i_1}),\gamma_1$ is a regular cardinal,
\sn
\item[${{}}$]  $(\beta) \quad \gamma_2$ has cofinality $>
  \kappa^\theta$ and $\gamma < \gamma_2 \Rightarrow
  |\gamma|^{\kappa^\theta} < \gamma_1$.
\end{enumerate}
\end{theorem}

\begin{remark}
\label{z37}
1) This relates on the one hand to the 
partition theorem of \cite[Ch.XI]{Sh:f} continuing
Rubin-Shelah \cite{RuSh:117}, Shelah \cite[Ch.XI]{Sh:f} 
and on the other hand to Komjath-Shelah \cite{KoSh:796}; the latter is
continued in Gruenhut-Shelah \cite{GhSh:909} but presently this is not used.

\noindent
2) Now \ref{z35} is what we use but we can get a somewhat more general
 result - see \ref{z53}.

\noindent
3) In \ref{z35} the case $\gamma_1 = \gamma_2 > 
|\cT_1|$ is equivalent to $\gamma_1 = \gamma_2 = \infty$.
\end{remark}

\begin{PROOF}{\ref{z35}}
Let $\cC = \{\bar c:\bar c = \langle c_\varrho:\varrho \in
\proj_{\bold i_1}\rangle,c_{<>} = \bold c(\rt(\cT_1))$ and 
where $c_{\varrho} \in \Rang(\bold c)$ or
just $(\exists \eta \in \cT_1)(\varrho = \proj_{\bold i_2}(\eta)
\wedge c_\varrho = \bold c(\eta))\}$.  For transparency \wilog \, we
assume $\Rang(\bold c \rest \max(\cT_{\bold i_1})),\Rang(\bold c
\rest (\cT_{\bold i_1} \backslash \max(\cT_{\bold i_1}))$ are
disjoint.  Clearly $|\cC| \le \kappa^{|\proj(\bold i_1)|} =
\kappa^\theta < \lambda$.

Fix for a while $\bar c \in \cC$, first let
$\cT_{\bar c} = \{\eta \in \cT_1$: if $\nu \trianglelefteq \eta$ then $\bold
c(\nu) = \bold c_{\proj(\nu,\bold i_1)}\}$ so a subtree of $\cT_1$, i.e. a
downward closed subset noting that $\rt_{\cT_1} \in \cT_{\bar c}$.  

Second, for $\eta \in \cT_1$, let $X^1_{\bar c,\eta}$ 
be $\suc_{\cT_{\bar c}}(\eta)$ if $\eta \in \cT_{\bar c} \cap
\Dom(\bar{\bold I}_1)$ and this set is $\in \bold I_{1,\eta}$ and be
$\emptyset$ otherwise.  Let $\cT'_{\bar c} = \{\eta \in \cT_{\bar c}$:
if $\ell < \ell g(\eta)$ and $\eta \rest \ell \in \Dom(\bold I_1)$
then $\eta \rest (\ell +1) \notin X^1_{\bar c,\eta}$,
i.e. $\suc_{\cT_{\bar c}}(\eta \rest \ell) := \{\nu \in \suc_{\cT_1}
(\eta):\nu \in \cT_{\bar c}\} \ne \emptyset \mod \bold
I_{1,\eta}\}$, again $\cT'_{\bar c}$ is a subtree of $\cT_{\bar c}$,
moreover $\bold i_{2,\bar c} = (\cT'_{\bar c},\bar{\bold I} \rest \cT'_{\bar
  c})$ is a $\pit$.

Third, for $\eta \in \cT'_{\bar c},\Dp_{\bold i_1}(\eta) \in \Ord \cup
\{\infty\}$ is well defined and, now for $\eta \in \cT_1$, let
$X^2_{\bar c,\eta}$ be $\{\nu \in \suc_{\cT'_{\bar c}}(\eta):
\Dp_{\bold i_{2,\bar c}}(\nu) \ge \Dp_{\bold i_{2,\bar c}}(\eta)\} 
= \emptyset \mod \bold I_{1,\eta}\}$ if $\eta \in
\cT'_{\bar c} \cap \Dom(\bar{\bold I}_1),\Dp_{\bold i_{2,\bar c}}(\eta) 
< \infty$ and be $\emptyset$ otherwise.

If for some $\bar c \in \cC,\Dp_{\bold i_{2,\bar c}}
(\rt_{\cT'_{\bar c}}) \ge \gamma_2$ easily we are done, so toward a
contradiction assume this is not the case, so recalling $\cf(\gamma_2)
> |\cC|$ clearly $\gamma_* = \sup\{\Dp_{\bold i_{2,\bar c}}
(\rt_{\cT'_{\bar c}}) +1: \bar c \in \cC\} < \gamma_2$.  Now
for each $\eta \in \Dom(\bar{\bold I}_1)$ clearly all $X^1_{\bar
  c,\eta},X^2_{\bar c,\eta}$ are from $\bold I_{1,\eta}$ and their
number is $\le 2|\cC| < \lambda$ hence $X_\eta := \cup\{X^1_{\bar
  c,\eta} \cup X^2_{\bar c,\eta}:\bar c \in \cC\}$ belong to 
$\bold I_{1,\eta}$.

Hence $\bold i_3$ is an $\pit$ and $\bold i_1 \le \bold i_3$
where $\bold i_3 = \bold i(3) := \bold i_1 \rest \{\eta \in 
\cT_1$: if $\ell < \ell
g(\eta)$ and $\eta \rest \ell \in \Dom(\bar{\bold I}_1)$ then $\eta
\rest (\ell +1) \notin X_\eta\}$; moreover by the definition of
$\Dp_{\bold i_3}$ and the choice of $\bold i_3$, clearly
\mn
\begin{enumerate}
\item[$(*)_1$]   $(a) \quad \bold i_3$ is a pit; moreover $\bold i_1 
\le_{\pr} \bold i_3$ hence
\sn
\item[${{}}$]  $(b) \quad \eta \in \cT_{\bold i_3} \Rightarrow
\Dp_{\bold i_3}(\eta) = \Dp_{\bold i_1}(\eta)$.
\end{enumerate}
\mn
Define $h$ by
\mn
\begin{enumerate}
\item[$(*)_2$]  $h$ is a function from $\cT_{\bold i_1} \times \cC$
defined by
\sn
\begin{enumerate}
\item[$\bullet$]  $h(\eta,\bar c)$ is $-1$ if $\eta \in \cT_{\bold
i_1} \backslash \cT'_{\bar c}$
\sn
\item[$\bullet$]  $h(\eta,\bar c)$ is $\Dp_{\bold i_{2,\bar c}}(\eta)$
if $\eta \in \cT'_{\bar c}$ and $\Dp_{\bold i_{2,\bar c}}(\eta) < \gamma_*$
\sn
\item[$\bullet$]  $\Dp(\eta,\bar c) = \gamma_*$ if none of the above.
\end{enumerate}
\end{enumerate}
\mn
We now choose $(\bold c_n,h_n,\cX_n,\bar{\cY}_n,\cS_n)$ 
by induction on $n$ such that:
\mn
\begin{enumerate}
\item[$\boxplus$]  $(a)(\alpha) \quad \cX_n$ is a subset of
$\cup\{\proj_{\bold i_1}(m):m \le n\}$
\sn
\item[${{}}$]  $\,\,\,\,\,\,(\beta) \quad$ if $n=k+1$ then $\cX_k = \cX_n
\cap (\cup\{\proj_{\bold i_1}(m):m \le k\})$
\sn
\item[${{}}$]  $\,\,\,\,\,\,(\gamma) \quad \cS_n \subseteq \cX_n$
\sn
\item[${{}}$]  $(b)(\alpha) \quad h_n$ is a function with
domain $\cX_n \times \cC$ to $\gamma_* +1$
\sn
\item[${{}}$]  $\,\,\,\,\,\,(\beta) \quad \bold c_n$ is a function from
$\cX_n$ to $\Rang(\bold c)$
\sn
\item[${{}}$]  $(c) \quad \bar{\cY}_n = \langle \cY_{n,\gamma}:\gamma <
  \gamma_1\rangle$
\sn
\item[${{}}$]  $(d)(\alpha) \quad \cY_{n,\gamma}$ is a subset of
$\cT_{\bold i_3}$, downward closed of cardinality $\le \theta$
\sn
\item[${{}}$]  $\,\,\,\,\,\,(\beta) \quad$ if $\eta \in \cY_{n,\gamma}$ then
  $\ell g(\eta) \le n$
\sn
\item[${{}}$]  $\,\,\,\,\,\,(\gamma) \quad$ if $\eta \in \cY_{n,\gamma}$ then
$\Dp_{\bold i_3}(\eta) = \Dp_{\bold i_1}(\eta) \ge \gamma$
\sn
\item[${{}}$]  $\,\,\,\,\,\,(\delta) \quad$ if $\eta \in \cY_{n,\gamma}$ and
$\ell g(\eta) < n$ and $\eta \notin \Dom(\bar{\bold I}_1)$ then

\hskip30pt $\suc_{\cT_{\bold i_3}}(\eta) = \suc_{\cT_{\bold i_1}}(\eta)$ is
  $\subseteq \cY_{n,\gamma}$
\sn
\item[${{}}$]  $\,\,\,\,\,\,(\varepsilon) \quad$ if $\eta \in 
\cY_{n,\gamma}$ and $\ell g(\eta) < n$ 
and $\eta \in \Dom(\bar{\bold I}_1)$ \then \, $\suc_{\cT_{\bold i_3}}(\eta)$

\hskip30pt  is a singleton 
\sn
\item[${{}}$]  $\,\,\,\,\,\,(\zeta) \quad$ if $\gamma < \gamma_2$ \then \,
$\cX_n = \{\proj_{\bold i_1}(\eta):\eta \in \cY_{n,\gamma}\}$
\sn
\item[${{}}$]  $\,\,\,\,\,\,(\eta) \quad$ if $\eta \in \cY_{n,\gamma}$ and
  $\nu = \proj_{\bold i _3}(\eta)$ \then \,:
\sn
\begin{enumerate}
\item[${{}}$]  $\bullet_1 \quad \bold c(\eta) = \bold c_n(\nu)$
\sn
\item[${{}}$]  $\bullet_2 \quad h_n(\nu,\bar c) = h(\eta,\bar c)$ for
  every $\bar c \in \cC$
\sn
\item[${{}}$]  $\bullet_3 \quad \eta \in \Dom(\bold I_1)$ \Iff \, $\nu
  \in \cS_n$.
\end{enumerate}
\sn
\item[${{}}$]  $\,\,\,(\theta) \quad$ follows: the function $\eta \mapsto
  \proj_{\bold i_3}(\eta)$ on $\cY_{n,\gamma}$ is one to one.
\end{enumerate}
\mn
\underline{Why this is possible}:

For $n=0$ this is trivial.

For $n=m+1$ for every $\gamma < \gamma_1$, choose $\bar\varrho_{n,\gamma} \in
\Pi\{\suc_{\cT_{\bold i_3}}(\eta):\eta \in \cY_{n,\gamma +1},\ell g(\eta)
= m,\eta \in \Dom(\bar{\bold I}_1)\}$ such that if $\eta \in
\Dom(\bar\varrho_{n,\gamma})$ then $\dep_{\bold i_1}(\eta) \ge \gamma$,
possible as $\eta \in  \Dom(\bar\varrho_{n,\gamma}) \Rightarrow
\dep_{\bold i_1}(\eta) \ge \gamma +1$.  Let $\cY'_{n,\gamma} =
\cY_{m,\gamma +1} \cup\{\nu$: for some $\eta \in \cY_{n,\gamma +1}$ we
have $\ell g(\eta) = m$ and we have $\eta \notin \Dom(\bar{\bold I}_1)
\Rightarrow \nu = \varrho_{n,\gamma}(\eta)$ and 
$\eta \notin \Dom(\bar{\bold I}_1) \Rightarrow \nu \in 
\suc_{\cT_{\bold i_3}}(\eta)\} \cup \Rang(\bar\varrho_{n,\gamma})$. 

Let $\cX'_{n,\gamma} = \{\proj_{\bold i_1}(\eta):\eta \in
\cY'_{n,\gamma}\}$ and let the function $\bold c'_{n,\gamma}:\cX'_{n,\gamma}
\rightarrow \Rang(\bold c)$ be defined by $\eta \in 
\cY'_{n,\gamma} \Rightarrow \bold c'_{n,\gamma}(\proj_{\bold
i_3}(\eta))) = \bold c(\eta)$, well defined as in $\boxplus(d)(\eta)$
and let $\cS_{n,\gamma} = \{\proj_{\bold i_1}(\eta):\eta \in
\cY'_{n,\gamma}$ and $\eta \in \Dom(\bar{\bold I}_1)\}$.   
Let $h_{n,\gamma}:\cX'_{n,\gamma} 
\rightarrow \gamma_* +1$ be defined by :
if $\bar c \in \cC,\nu = \proj_{\bold i_{2,\bar c}}(\eta)$ and $\eta \in
\cY_{n,\gamma}$ then $\eta \notin \cT'_{\bar c} \Rightarrow
h_{n,\gamma}(\nu) = \gamma,\eta \in \cT'_{\bar c} \Rightarrow
h_{n,\gamma}(v) = \Dp_{\bold i_{2,\bar c}}(\eta)$.

Now $\cX'_{n,\gamma}$ is a subset of $\proj_{\bold i_1}$, a set of
cardinality $\le \theta$ and $\bold c'_{n,\gamma}$ is a funtion from
$\cX'_{n,\gamma}$ into $\Rang(\bold c)$, a set of cardinality $\le
\kappa$ and $h_{n,\gamma}$ is a function from $\cX'_{n,\gamma}
\subseteq \proj_{\bold i_1}$ into $\gamma_*$.  But $\gamma_* <
\gamma_2,\gamma_* + \kappa < \gamma_1,\gamma_1$ is a regular cardinal
(recalling clause (e) of the theorem) 
and $(|\gamma_*| + \kappa)^\theta < \cf(\gamma_1) = \gamma_1$ hence
for every $\gamma < \gamma_1$ we have $|\{(X'_{n,\gamma},\bold
c_{n,\gamma},h_{n,\gamma}):\gamma < \gamma_1\}| \le 2^\theta \cdot
\kappa^\theta \cdot |\gamma_*|^\theta < \cf(\gamma_1) = \gamma_1$
hence for some $\bold c_n,h_n,\cX_n$ the set $S_n := \{\gamma <
\gamma_1:\bold c'_{n,\gamma} = \bold c_n$ and $h_{n,\gamma} =
h_n,\cX'_{n,\gamma} = \cX_n$ and
$\cS_{n,\gamma} = \cS_n\}$ is unbounded in $\gamma_1$.

Lastly, let $\cY_{n,\gamma} = \cY'_{n,\min(S_n \backslash \gamma)}$,
clearly $\bold c_{n+1},h_{n+1},\langle \cY_{n,\gamma}:\gamma <
\gamma_2\rangle$ are as required; so we can carry the induction.
\smallskip

\noindent
\underline{Why this is enough}:

Let $\cX = \cup\{\cX_n:n < \omega\} \subseteq \proj(\bold i_1)$ and
$\cS = \cup\{\cS_n:n < \omega\}$ and 
$\bold c = \cup\{\bold c_n:n < \omega\}$ and $\bold h = \cup\{h_n:n <
\omega\}$ so by $\boxplus(d)(\eta)$ clearly there is $\bar c_* \in \cC$ such
that $c_\varrho = c(\varrho)$ when the latter is defined, so: 
\mn
\begin{enumerate}
\item[$\odot_1$]  if $n < \omega,\gamma < \gamma_1,\eta \in
\cY_{n,\gamma}$ and $\nu = \proj(\bold i_1) \in \cX$ \ then \,
\sn
\begin{enumerate}
\item[$(a)$]   $\bold c(\eta) = \bold c_n(\proj_{\bold i_1}(\eta))$
\sn
\item[$(b)$]  $\Dp_{\bold i_{2,\bar c}}(\eta) = h(\eta,\bar c) =
  h_n(\nu,\bar c)$
\sn
\item[$(c)$]  $\Dp_{\bold i_1}(\eta) \ge \gamma$
\end{enumerate}
\end{enumerate}
\mn
Also
\mn
\begin{enumerate}
\item[$\odot_2$]  $\cX \subseteq \proj_{\bold i_1}$ is a set of 
finite sequences, closed under initial segments 
with no $\triangleleft$-maximal member.
\end{enumerate}
\mn
[Why?  Straight, e.g. if $\nu \in X$ choose $n = \ell g(\nu)+2$ let
  $\gamma < \gamma_1$ and choose $\eta \in Y_{n,\gamma +1}$ such that
  $\proj_{\bold i_1}(\eta)=\nu$, now by clause (c) of $\odot_1$ we
  know that $\Dp_{\bold i_1}(\eta) \ge \gamma +1$, hence there is
  $\eta_1 \in \suc_{\cT_{\bold i_1}}(\eta)$ in $Y_{n,\gamma +1}$ hence
  $\nu_1 = \proj_{\bold i_1}(\eta_1)$ is in $\suc_{\cX}(\nu)$,
  i.e. successor of $\eta$ in $\cX_{n+1}$ hence in $\cX$.]
\mn
\begin{enumerate}
\item[$\odot_3$]   if $\nu \in \cX$ then $\bold h(\nu,\bar c) \ne -1$.
\end{enumerate}
\mn
[Why?  Let $n > \ell g(\nu)$, let $\gamma < \gamma_2$.  Now by
$\boxplus(d)(\zeta)$  there is $\eta \in \cY_{n,\gamma}$ such
that $\proj_{\bold i_1}(\eta)=\nu$.

Next by $(*)_2$ we have $h(\eta,\bar c)$ is -1 iff $\eta \notin
\cT'_{\bar c}$.  However, $\eta \in \cT_{\bar c}$ by the definition of
$\cT_{\bar c}$ and the choice of $\bar c$ and $\boxplus(d)(\eta)$;
moreover $\eta \in \cT'_{\bar c}$ by the definition of $\cT'_{\bar c}$ anda
of $\bold i_3$ and clause $\boxplus(d)(\alpha)$.

By the last two sentences $h(\eta,\bar c) \ne -1$ hence by the choice
of $\eta$, i.e. as $\proj_{\bold i_1}(\eta)=\nu$, clause
$\boxplus(d)(\eta)$ tells us $\bold h(\nu,\bar c) = h(\eta,\bar c)$ so
together $\bold h(\nu,\bar c) \ne -1$ as promised.]
\mn
\begin{enumerate}
\item[$\odot_4$]  $0 \le \Dp_{\bold i_{2,\bar c}}(\langle \rangle) <
  \gamma_*$ hence $\bold h(<>,\bar c) < \gamma_*$.
\end{enumerate}
\mn
[Why?  Similarly using $\boxplus(d)(\eta)\bullet_3$.]
\mn
\begin{enumerate}
\item[$\odot_5$]  if $\nu \in \cX \backslash \cS$ and $0 \le \bold
h(\nu,\bar c) < \gamma_*$ \then \, for some $\rho \in \suc_{\cX}(\nu)$
we have $0 \le \bold h(\rho,\bar c) < \bold h(\nu,\bar c) < \gamma_*$.
\end{enumerate}
\mn
[Why?  Similarly using $\boxplus(d)(\delta)$.]
\mn
\begin{enumerate}
\item[$\odot_6$]  if $\nu \in \cS(\subseteq \cX)$ and $0 \le \bold
h(\nu,\bar c) < \gamma_*$ \then \, for the unique $\rho \in
\suc_{\cX}(\nu)$ we have $0 \le \bold h(\rho,\bar c) < \bold
h(\nu,\bar c) < \gamma_*$.
\end{enumerate}
\mn
[Why?  Similarly using $\boxplus(d)(\varepsilon)$.]

By $\odot_4,\odot_5,\odot_6$ together we get a contradiction.
\end{PROOF}

\noindent
We may prefer the following variant of \ref{z35}.
\begin{definition}
\label{z51}
1) For a pit $\bold i = (\cT,\bar{\bold I})$ and partition $\bar{\cS}
   = (\cS_0,\cS_1)$ of $\cS_{\bold i}$ (or just $\bar{\cS} =
   (\cS_0,\cS_1)$ such that $\cS_0 \cap \cS_1 = \emptyset$ and
   $\cS_{\bold i} \subseteq \cS_0 \cup \cS_1$) we define $\Dp_{\bold
   i,\cS}:\cT \rightarrow \Ord \cup \{\infty\}$, stipulating $\infty
   +1=\infty$ by defining when $\Dp_{\bold i,\lambda}(\eta) \ge
   \alpha$ by induction on the ordinal $\alpha$ (compare with
   \ref{z32}(4)):
\mn
\begin{enumerate}
\item[$(a)$]  if $\eta \in \max(\cT)$ \then \, $\Dp_{\bold
  i,\lambda}(\eta) \ge \alpha$ iff $\alpha = 0$
\sn
\item[$(b)_0$]  if $\eta \in \cS_0$ hence $\eta \in \cS,\eta \notin 
\max(\cT)$ \then \, $\Dp_{\bold i,\bar{\cS}}(\eta) \ge \alpha$ iff for
every $\beta < \alpha$ the set $\{\nu \in \suc_{\cT}(\eta):\Dp_{\bold
  i,\bar{\cS}}(\nu) \ge \beta\}$ belong to $\bold I^+_\eta$
\sn
\item[$(b)_1$]  if $\eta \in \cS_1$ \then \, $\Dp_{\bold
  i,\bar{\cS}}(\eta) \ge \alpha$ iff $\{\nu \in
  \suc_{\cT}(\eta):\Dp_{\bold i,\bar{\cS}}(\nu) \ge \alpha\}$ belongs
  to $\bold I^+_\eta$
\sn
\item[$(c)$]  if $\eta \in \cT \backslash \cS \backslash \max(\cT)$ 
\then \, $\Dp_{\bold i,\bar{\cS}}(\eta) \ge \alpha$ iff for every $\nu \in
\suc_{\cT}(\eta)$ we have $\Dp_{\bold i,\bar{\cS}}(\eta)$.
\end{enumerate}
\end{definition}

\begin{theorem}
\label{z53}
There are a pit $\bold i_2$ and $\bar c = \langle c_\eta:\eta \in
 \proj(\bold i_1)\rangle$ such that $\bold i_1 \le \bold
 i_2,\Dp_{\bold i_2,\bar{\cS}}(\rt_{\bold i_2}) \ge \gamma_2$ and
 $\eta \in \cT_{\bold i_2} \Rightarrow \bold c(\eta) = \bold
 c_{\proj(\eta,\bold i_1)}$ \when \,:
\mn
\begin{enumerate}
\item[$(a)-(e)$]  as in \ref{z35} replacing $\Dp_{\bold i_2}$ by
  $\Dp_{\bold i_2,\bar{\cS}}$ in $(e)(\alpha)$
\sn
\item[$(f)$]  $\bar{\cS} = (\cS_0,\cS_1)$ is a partition of
  $\cS_{\bold i}$.
\end{enumerate}
\end{theorem}

\begin{PROOF}{\ref{z53}}
Similarly.
\end{PROOF}
\newpage

\section {Approximation to $\EM$ models} \label{Approx}

In the game below the protagonist tries to exemplify in a weak form
that the standard $\EM_{\tau(\gk)}(\lambda,\Phi)$ is
$\le_{\gk}$-embeddable into $N$ over $M$.  We may consider games in
which the protagonist tries to exemplify a weak form of isomorphism,
this is connected to logics which have $\EM$ models, continuing \cite{Sh:797},
but not for now.

Here we do not try to get the best cardinal bounds; just
enough for the result promised in the abstract.

\begin{definition}
\label{a2}
Assume $\lambda > \kappa \ge \LST_{\gk} + |\tau_{\gk}|$ and $M \in
K^{\gk}_\kappa$ and $M \le_{\gk} N$ and $\gamma$ is an ordinal.

\noindent
1) We say $\Phi$ is an $(M,\lambda,\kappa,\gamma)$-solution of $N$ or is an
$(N,M,\lambda,\kappa,\gamma)$-solution \when \, $\Phi \in
\Upsilon^{\sor}_\kappa(\gk_M)$ and in the game
$\Game^1_{N,M,\lambda,\Phi,\gamma}$ the protagonist has a winning strategy.

\noindent
2) Assume $\Phi \in \Upsilon_\kappa(\gk_M)$ recalling Definition
   \ref{z12} fixing $M_\lambda =
   \EM(\lambda,\Phi)$ and $M_I = \EM(I,\Phi)$ for $I \subseteq
   \lambda$ and \wilog \, every $M_I$ (equivalently some $M_I$) 
is standard, hence in particular 
$M \le_{\gk} M_I \rest \tau(\gk)$.  We define the game
$\Game^1_{N,M,\lambda,\Phi,\gamma}$, a play last $< \omega$ moves, in the
$n$-th move $\lambda_n,J_n,\bar h_n,\gamma_n$ are chosen such that:
\mn
\begin{enumerate}
\item[$\boxplus_n$]  $(a) \quad \lambda_0 = \lambda$
\sn
\item[${{}}$]  $(b) \quad$ if $n=m+1$ then $\kappa < \lambda_n <
  \lambda_m$ moreover $\lambda_m \rightarrow (\lambda_n)^n_{2^\kappa}$
\sn
\item[${{}}$]  $(c) \quad J_0 = \lambda$, and if $n=m+1$ then $J_n
  \subseteq J_m$
\sn
\item[${{}}$]  $(d) \quad |J_n| = \lambda_n$
\sn
\item[${{}}$]  $(e) \quad \bar h_n = \langle h_u:u \in [J_n]^n\rangle$
\sn
\item[${{}}$]  $(f) \quad$ if $u \in [J_n]^n$ \then \, $h_u$ is a
  $\le_{\gk}$-embedding of $M_u$ into $N$ extending $h_v$ 

\hskip25pt whenever $v \subseteq u$

\begin{scriptsize}
[Explanation: note if $v \subset u,|v| = m$ then $v \in [J_n]^m
  \subseteq [J_m]^n$ hence $h_v$ was defined; this says then for
  $u_1,u_2 \in [J_n]^n,h_{u_2},h_{u_2}$ are compatible functions]
\end{scriptsize}
\sn
\item[${{}}$]  $(g) \quad \gamma_0 = \gamma$ and $\gamma_{n+1}$ is an
  ordinal $< \gamma_n$.
\end{enumerate}
\mn
In the $n$-th move:
\mn
\begin{enumerate}
\item[$(A)$]  if $n=0$ the antagonist chooses $\lambda_0 = \lambda,J_0
 = \lambda,\gamma_0 = \gamma$ and the protagonist chooses $\bar h_0$
\sn
\item[$(B)$]  if $n=m+1$ then
\sn
\begin{enumerate}
\item[$(a)$]  the antagonist chooses an ordinal $\gamma_n < \gamma_m$
 and $\lambda_n >\kappa$ such that
  $\lambda_m \rightarrow (\lambda_n)^n_{\beth_2(\kappa)}$
\sn
\item[$(b)$]  the protagonist chooses $\bar h_n$ and $\cS_n \in
(\ER^n_{J_m,\lambda_n,\beth_2(\kappa)})^+$, i.e. $\cS_n \subseteq
[\lambda_m]^{\lambda_n}$ and $\cS_n$ is not from this ideal,
see Definition \ref{z18}
\sn
\item[$(c)$]   the antagonist chooses $J_n \in \cS_n \subseteq 
[J_m]^{\lambda_n}$
\end{enumerate}
\sn
\item[$(C)$]  the play ends when a player has no legal move and then
this player loses.
\end{enumerate}
\end{definition}

\noindent
Another presentation:
\begin{definition}
\label{a5}
Assume $M \le_{\gk} N$ and $\LST_{\gk} + |\tau_{\gk}| \le \theta,\|M\|
+ \theta \le \kappa < \lambda$ and $\Phi \in 
\Upsilon^{\oor}_\theta[M,\gk]$.

\noindent
1) Below we omit $\gamma$ if (a) or (b), where:
\mn
\begin{enumerate}
\item[$(a)$]  $\gamma = \cf(\lambda),\lambda$ strong limit and $\alpha
< \cf(\lambda) \Rightarrow |\alpha|^{2^{\kappa + \|M\|}} <
\cf(\lambda)$
\sn
\item[$(b)$]  not (a) but $\gamma$ is maximal
such that $\gamma = \omega \gamma$ is infinite 
and $\beth_\gamma(\kappa + \|M\|) \le \lambda$ and $\lambda$ is 
strong limit of cofinality $> \beth_2(\kappa)$ 
\newline
(similarly in all such definitions).
\end{enumerate}
\mn
2) We say that $\bold x$ is a direct witness for
$(N,M,\lambda,\kappa,\gamma,\Phi)$ \when \, $\bold x$ consists of:
\mn
\begin{enumerate}
\item[$(a)$]  $N,M,\Phi,\lambda,\kappa$ and $\gamma$
\sn
\item[$(b)$]  $\cT$ is a non-empty set of finite sequences closed
  under initial segments
\sn
\item[$(c)$]  if $\eta \in \cT$ \then \,:
\sn
\item[${{}}$]  $(\alpha) \quad \eta(2n)$ is a cardinal when $2n < \ell g(\eta)$
\sn
\item[${{}}$]  $(\beta) \quad 
\eta(2n+1)$ is a subset of $\lambda$ of cardinality
  $\eta(2n)$ when $2n+1 < \ell g(\eta)$
\sn
\item[${{}}$]  $(\gamma) \quad 
\eta(2n+1) \supseteq \eta(2n+3)$ when $2n+3 < \ell g(\eta)$
\sn
\item[${{}}$]  $(\delta) \quad \eta(2n) \ge \eta(2n+2)$,
moreover $\eta(2n) \rightarrow (\eta(2n+2))^{2n+1}_{\beth_2(\kappa)}$ when

\hskip30pt  $2n+2 < \ell g(\eta)$
\sn
\item[$(d)$]  $I_\eta,\lambda_\eta$ for $\eta \in \cT$ are defined by:
\sn
\begin{enumerate}
\item[$(\alpha)$]  if $\ell g(\eta) = 0$ then 
$I_\eta = \lambda,\lambda_\eta = \lambda$
\sn
\item[$(\beta)$]  if $\ell g(\eta) = 2n+1$ then $I_\eta = I_{\eta
  \rest (2n)}$, see $(\alpha)$ or $(\gamma)$ and $\lambda_\eta = \eta(2n)$
\sn
\item[$(\gamma)$]  if $\ell g(\eta) = 2n+2$ then $I_\eta =
 \eta(2n+1),\lambda_\eta = \eta(2n) = \lambda_{\eta \rest (2n+1)}$,
 see $(\alpha)$ or $(\beta)$
\end{enumerate}
\sn
\item[$(e)$]  if $\eta \in \cT \backslash \max(\cT)$ has length $2n+1$ then:
the set $\cS_\eta = \{I_\nu:\nu \in \suc_{\cT}(\eta)\} 
\subseteq [I_\eta]^{\lambda_\eta}$ is not from the ideal
  $\ER^{\lfloor \ell g(\eta)/2\rfloor}_{I_\eta,\lambda_\eta,\beth_2(\kappa)}$
\sn
\item[$(f)$]  if $\eta \in \cT$ then:
\sn
\begin{enumerate}
\item[$(\alpha)$]  $\bar h_\eta = \langle h_{\eta,u}:u \in
  [I_\eta]^{\le \lfloor \ell g(\eta)/2\rfloor}\rangle$
\sn
\item[$(\beta)$]  $h_{\eta,u}$ is a $\le_{\gk}$-embedding of 
$\EM_{\tau(\gk)}(u,\Phi)$ into $N$ for $u \in [I_\eta]^{\le \lfloor 
\ell g(\eta)/2 \rfloor}$
\sn
\item[$(\gamma)$]  $u_1 \subseteq u_2 \in [I_\eta]^{\le \lfloor \ell
g(\eta)/2\rfloor} \Rightarrow h_{\eta,u_1} \subseteq h_{\eta,u_2}$
\sn
\item[$(\delta)$]  if $u \in [I_\eta]^{\lfloor \ell g(\eta)/2\rfloor}$ and $\nu
 \triangleleft \eta$ and $\ell g(\nu) \ge 2|u|$, \then \,
  $h_{\eta,u} = h_{\nu,u}$
\sn
\item[$(\varepsilon)$]  if $\ell g(\eta) = 2n+2$ and $u \in
[I_\eta]^{\le n}$ then $h_{\eta,u} = h_{\eta \rest (2n+1),u}$
\sn
\item[$(\zeta)$]  there is $\bar a = \bar a_{\bold x} = \langle
  a_\alpha:\alpha < \lambda\rangle \in {}^\lambda N$ such that $\alpha
\in u \in [I_\eta]^{\le \lfloor \ell g(\eta)\rfloor/2} +
h_{\eta,u}(\alpha) = a_\alpha$ and $\bar a$ is with no repetitions
\end{enumerate}
\sn
\item[$(g)$]  $\Dp_{\bold x}(<>) \ge \gamma$ where $\Dp_{\bold
  x}(\eta)$ is defined as $\Dp_{\bold i(\bold x)}(\eta)$, see
Definition \ref{z32}, where $\bold i = \bold i(\bold x) = \bold
  i_{\bold x}$ is defined by:
\sn
\begin{enumerate}
\item[$\bullet$]  $\cT_{\bold i} = \cT$
\sn
\item[$\bullet$]  $\cS_{\bold i} = \{\eta \in \cT:\eta$ is not
$\triangleleft$-maximal in $\cT$ and $\ell g(\eta)$ is $\odd\}$
\sn
\item[$\bullet$]  if $\eta \in \cS_{\bold i}$ and $\ell g(\eta)$ is
  odd then $\bold I_{\bold i,\eta} = 
\ER^{\lfloor \ell g(\eta)\rfloor}_{I_\eta,\lambda_\eta,\beth_2(\kappa)}$
recalling \ref{z32}(1A)
\sn
\item[$\bullet$]  if $\eta \in \cS_{\bold i}$ and $\ell g(\eta)$ is even
then $\cI_{\bold i,\eta} = \{\emptyset\}$.
\end{enumerate}
\end{enumerate}
\mn
\end{definition}

\begin{definition}
\label{a6}
1) We say $\bold x$ is a pre-$\gk$-witness of 
$(N,M,\lambda,\kappa,\delta)$ \when \, it as in \ref{a5} omitting
$\bar h$, i.e. clause (f), so $N,M$ are irrelevant.

\noindent
2)  We say $\bold x$ is a semi-$\gk$-witness of 
$(N^+,M,\lambda,\kappa,\delta)$ \when \,: it consists of:
\mn
\begin{enumerate}
\item[$(a)$]  $N^+$ expands a model from $K_{\gk},M \le_{\gk}(N^+
  \rest \tau(\gk)),\lambda \ge \kappa \ge (\tau(N^+))$
\sn
\item[$(b)-(e)$]  as in \ref{a5}(2)
\sn
\item[$(f)$]  $\bar a = \langle a_\alpha:\alpha < \lambda\rangle$
\sn
\item[$(g)$]  as in \ref{a5}(2).
\end{enumerate}
\end{definition}

\begin{claim}
\label{a7}
1) The definitions \ref{a2}, \ref{a5} are equivalent.

\noindent
2) In Definition \ref{a5}, $\bold i_{\bold x}$ is indeed a pit.

\noindent
3) If $\Phi_1 \bold E^{\ai}_\kappa \Phi_2,\Phi_\ell \in
   \Upsilon^{\sor}_\kappa[M,\gk]$ for $\ell=1,2$ and $\Phi_1$ has a
   $(N,M,\lambda,\kappa)$-witness \then \, $\Phi_2$ has a
   $(N,M,\lambda,\kappa)$-witness. 
\end{claim}

\begin{PROOF}{\ref{a5}}
Straightforward.
\end{PROOF}

\begin{claim}
\label{a8}
1) If $\Phi_\ell \in \Upsilon^{\sor}_\kappa[\gk_M],
\kappa \ge \tau(\gk) + \|M\|$
 and $M_\ell = \EM_{\tau(\gk)}(\lambda,\Phi_\ell)$ for 
$\ell=1,2$ and $\lambda$ is
strong limit of cofinality $\mu$ where $\mu = (\beth_2(\kappa))^+$ or $\mu$ is
 regular such that
$(\forall \alpha < \mu)(|\alpha|^{2^\kappa} < \mu)$ and the 
protagonist wins in the game $\Game^1_{M_2,M,\lambda,\Phi_1,\mu}$
(equivalently some $\bold x$ is a witness for
$(M_2,M,\lambda,\kappa,\Phi_1)$) \then \, $\Phi_1 \le^3_\kappa \Phi_2$,
see Definition \ref{z24}.
\end{claim}

\begin{PROOF}{\ref{a8}}
Straightforward by \ref{z35} and the definitions of the ideal $\ER$ in
\ref{z18}.  See details in a similar case in the proof of \ref{a9}(1) below.
\end{PROOF}

\begin{claim}
\label{a9}
Assume $M \le_{\gk} N,\kappa \ge \|M\| + \theta,\theta \ge
\LST_{\gk} + |\tau_{\gk}|$
and $\|N\| \ge \lambda,\lambda$ strong limit of 
cofinality $\mu$ and $\mu = 
(\beth_2(\kappa))^+$ or $\mu$ is regular such that $(\forall
\alpha < \mu)(|\alpha|^{2^\kappa} < \mu)$.

\noindent
1) There are $\bold x,\Phi$ such that:
\mn
\begin{enumerate}
\item[$(a)$]  $\Phi \in \Upsilon^{\sor}_\theta(\gk_M)$
\sn
\item[$(b)$]  $\bold x$ is a direct witness of $(N,M,\lambda,\kappa,\Phi)$.
\end{enumerate}
\mn
2) If $M_1 = M,\Phi_1 \in \Upsilon^{\sor}_\theta[\gk_{M_1}]$ and $\bold x_1$
a direct witness for $(N,M_1,\lambda,\kappa,\Phi_1)$ and $M_1 \le_{\gk} M_2
  \le_{\gk} N$ and $\|M_2\| \le \kappa$ \then \, there are
  $\Phi_2,\bold x_2$ such that:
\mn
\begin{enumerate}
\item[$(a)$]  $\Phi_2 \in \Upsilon^{\sor}_\theta[M_2]$
\sn
\item[$(b)$]  $\Phi_1 \le^1s_\kappa \Phi_2$ and $\Phi_1 \le^4_\kappa \Phi_2$
\sn
\item[$(c)$]  $\bold x_2$ is a direct witness $(N,M_2,\lambda,\kappa,\Phi_2)$.
\end{enumerate}
\mn
3) If in part (1) we change the assumption on $\lambda$ to $\lambda =
\beth_{\omega \cdot \gamma}(\kappa)$ \then \, there are $\Phi,\bold x$
such that:
\mn
\begin{enumerate}
\item[$(a)$]  $\Phi \in \Upsilon^{\sor}_\kappa[M,\gk]$
\sn
\item[$(b)$]  $\bold x$ is a direct witness of
  $(N,M,\Phi,\lambda,\kappa,\gamma,\Phi)$. 
\end{enumerate}
\mn
4) Also part (2) has a version with $(\gamma_1,\gamma_2)$ as in \ref{z35}.
\end{claim}

\begin{PROOF}{\ref{a9}}
1) Let $\langle a_\alpha:\alpha < \lambda\rangle$ be a sequence of
   pairwise distinct members of $N$.

Now
\mn
\begin{enumerate}
\item[$(*)_1$]  let $\cT$ be the set of finite sequences $\eta$ satisfying
clauses (b),(c) of Definition \ref{a5}
\sn
\item[$(*)_2$]  let $\bar{\bold I} = \langle \bold I_\eta:\eta \in 
\cS \rangle$ where
\sn
\item[${{}}$]  $\bullet \quad \cS = \{\eta \in \cT:\eta$ is not
  $\triangleleft$-maximal in $\cT\}$
\sn
\item[${{}}$]  $\bullet \quad$ if $\eta \in \cS,\ell g(\eta) = 2n+1$
  then $\bold I_\eta = \ER^n_{I_\eta,\lambda_\eta,\beth_2(\kappa)}$
\sn
\item[${{}}$]  $\bullet \quad$ if $\eta \in \cS$ and $\ell g(\eta) = 2n$
then $\bold I_\eta = \{\emptyset\}$, the trivial ideal
\sn
\item[$(*)_3$]  $\bold i_1 = \bold i(1) = (\cT,\bar{\bold
  I})$ is a pit and is $(2^\kappa)^+$-complete and $\Dp_{\bold
  i_1}(<>) \ge (\beth_2(\kappa))^+$.
\end{enumerate}
\mn
[Why?  Just read Definition \ref{z32}(3) and the ideal $\ER$ is from
  Definition \ref{z18} and it is $(2^\kappa)^+$-complete by \ref{z8d}
  and as for the depth recall $\mu = (\beth_2(\kappa))^+$.]
\mn
\begin{enumerate}
\item[$(*)_4$]  Let $M^+$ be such that:
\sn
\begin{enumerate}
\item[$(a)$]  $M^+$ is an expansion of $N$
\sn
\item[$(b)$]  $|\tau(M^+)| \le \kappa$ and $\tau' := \tau(M^+) \backslash
  \{c_a:a \in M\}$ has cardinality $\le \theta$
\sn
\item[$(c)$]  if $M^+_1 \rest \tau' \subseteq M^+ \rest \tau'$ 
then $M^+_1 \rest \tau(\gk) \le M^+ \rest \tau(\gk)$
\sn
\item[$(d)$]  $|M| = \{c^{M^+}:c \in \tau(M^+)\}$.
\end{enumerate}
\end{enumerate}
\mn
[Why $M^+$ exists?  By the representation theorem, \cite[\S1]{Sh:88r}
except clause (d) which as before is easy.]

We like to apply Theorem \ref{z35} but before this we need
\mn
\begin{enumerate}
\item[$(*)_5$]  there is a pit $\bold i_2 = \bold i(2)$ such that
$\bold i(1) \le_{\pr} \bold i(2)$ (see \ref{z32}(2A)) so $\Dp_{\bold
  i(2)}(\eta) = \Dp_{\bold i(1)}(\eta)$ for $\eta \in \cT_{\bold i(2)}$ and:
\sn
\begin{enumerate}
\item[$\bullet$]   if $\eta \in \cT_{\bold i(2)},\ell g(\eta) = 2n+1$
  and $\nu \in \suc_{\cT_{\bold i(2)}}(\eta)$ then $\langle
  a_\alpha:\alpha \in \nu(2n+1)\rangle$ is an $n$-indiscernible
  sequence in $M^+$ for quantifier free formulas, may add: and $N \rest
  \{\sigma_\varepsilon(a_{\alpha_0},\dotsc,a_{\alpha_{n-1}}):\varepsilon
  < \zeta\} \le_{\gk} N$ when $\zeta < \kappa^+$ and
  $\sigma_\varepsilon$ is a $\tau(M^+)$-term.
\end{enumerate}
\end{enumerate}
\mn
[Why such $\bold i(2)$ exists?  By the definition of the ideal $\bold
  I_\eta$, see $(*)_2$ above and by Definition \ref{z8}.  
That is, for $\eta \in \Dom(\bold I_{\bold i_1})$ of 
length $2n+1$ let $X_\eta =\{\nu:\nu
  \in \suc_{\cT}(\eta),\langle a_\alpha:\alpha \in \nu(2n+1)\rangle$ is
$n$-indiscernible in $M^+$ for quantifier free formulas$\}$,
recalling $\Dom(\bold I_{\bold i_1,\eta}) = 
\{u \subseteq I_\eta:|u| = \eta(2n)\}$.  By \ref{z18} clearly $X_\eta
= [\lambda_\eta]^{\eta(2n)} \mod
\ER^n_{\lambda_\eta,\eta(2),\beth_2(\kappa)}$, see Definition \ref{z32}(1A).

Now let $\cT' = \{\eta \in \cT$: if $2n+1 <
\ell g(\eta)$ then $\eta \rest (2n+2) \in X_\eta\}$ and $\bold i_2 =
\bold i_1 \rest \cT'$, so clearly $\bold i_1 \le_{\pr} \bold i_2$, see
Definition \ref{z32}(2A).]

Next
\mn
\begin{enumerate}
\item[$(*)_6$]  define a function $\bold c$ with domain $\cT_{\bold
  i_2}$ as follows: 
\sn
\begin{enumerate}
\item[$\bullet$]  if $\eta \in \cT,\ell g(\eta) = 2n+2$, \then \,
  $\bold c(\eta)$ is the quantifier type in $M^+$ of $\langle
  a_\ell:\ell < n \rangle$ for any $\alpha_0 < \alpha_1 < \ldots <
\alpha_{n-1}$ from $\eta(2n+1)$
\sn
\item[$\bullet$]  if $\eta \in \cT,\ell g(\eta) = 2n+1$ or $\ell
g(\eta) = 0$, \then \, $\bold c(\eta) = 0$.
\end{enumerate}
\end{enumerate}
\mn
Clearly
\mn
\begin{enumerate}
\item[$(*)_7$]  $\Rang(\bold c)$ has cardinality $\le 2^\kappa = 2^\kappa$.
\end{enumerate}
\mn
So by \ref{z35} (with a degenerate projection; so $\kappa,\theta$
there stands for $2^\kappa,\aleph_0$ here):
\mn
\begin{enumerate}
\item[$(*)_8$]  there are $\bold i(3) = \bold i_3 \ge \bold i_2$ and
  $\langle c_n:n < \omega\rangle$ such that:
\sn
\item[${{}}$]  $(a) \quad \eta \in \cT_{\bold i_3} 
\Rightarrow \bold c(\eta) = c_{\ell g(\eta)}$
\sn
\item[${{}}$]  $(b) \quad \Dp_{\bold i_3}(<>) \ge \beth_2(\kappa)$.
\end{enumerate}
\mn 
The rest should be clear.

\noindent
2) Similar proof, this time in $M^+$ we have individual constants for
   every member of $M_2$ and we start with the witness $\bold x_1$ so
$X_\eta$ have fewer elements still positive modulo the ideal.

\noindent
3),4) Similarly.
\end{PROOF}

\begin{definition}
\label{a12}
We say $\bold x$ is an indirect witness for
$(N,M,\lambda,\kappa,\gamma,\Phi)$, recalling \ref{a5}(1), 
\when \, for some $\Psi$:
\mn
\begin{enumerate}
\item[$(a)$]  $M,N,\lambda,\kappa,\gamma,\Phi$ are as in Definition \ref{a5}
\sn
\item[$(b)$]  $\Psi \in \Upsilon^{\sor}_\kappa[\gk_M]$ and $\Phi \le^4_\kappa
  \Psi$, see Definition \ref{z24} 
\sn
\item[$(c)$]  $\bold x$ is a direct witness of
 $(N,M,\lambda,\kappa,\gamma,\Psi)$.
\end{enumerate}
\end{definition}

\begin{remark}
Why do we need the indirect witnesses?  As if we use direct witness only
in the proof of \ref{a22} it is not clear how to get many
non-isomorphic models.
\end{remark}

\begin{claim}
\label{a13}
If $I=I_\chi$ is as in \ref{z9}.

If (A) then (B) where:
\mn
\begin{enumerate}
\item[$(A)$]  $(a) \quad \LST_{\gk} + |\tau_{\gk}| \le 
\kappa < \chi_1 < \chi_2 < \chi_3 \le \chi$
  and for $\ell=\eta,2,\chi_{\ell+1}$ 

\hskip25pt is strong limit of cofinality $> \beth_2(\chi_\ell)$
\sn
\item[${{}}$]  $(b) \quad N = \EM_{\tau(\gk)}(I,\Phi_1)$ where $\Phi_1
  \in \Upsilon^{\sor}_\kappa[M_1,\gk],\|M_1\| \le \chi_1$
\sn
\item[${{}}$]  $(c) \quad M_2 \le_{\gk} N$ and $\|M_2\| \le \chi_1$
\sn
\item[${{}}$]  $(d) \quad \Phi_2 \in \Upsilon^{\sor}_\kappa[M,\gk]$
\sn
\item[${{}}$]  $(e) \quad \Phi_2$ has a witness for
 $(N,M_2,\chi_2,\kappa)$
\sn
\item[$(B)$]  $(a) \quad \Phi_2$ has a witness for
  $(N,M_2,\chi_3,\kappa)$
\sn
\item[${{}}$]  $(b) \quad$ if in addition $M_2 \le_{\gk} M_1$ then 
$\Phi_2 \le^3_\kappa \Phi_1$
\sn
\item[${{}}$]  $(c) \quad$ we can $\le_{\gk}$-embed
 $\EM_{\tau(\gk)}(I_\chi,\Phi_2)$ into $N$.
\end{enumerate}
\end{claim}

\begin{PROOF}{\ref{a13}}
As in the proof of \ref{a9} recalling the choice of $I$ in \ref{z9};
for (B)(c) we use Clause $(B)^+$ of \ref{a9}.
\end{PROOF}

\begin{remark}
In fact, in \ref{a13}, $\chi_2 = \beth_{1,1}(\chi_1)$ and $\chi_3 =
\beth_{\omega \gamma}(\chi_1)$ suffices so, of course, in (B)(a) we
use $(N,M_1,\chi_3,\kappa,\gamma)$.
\end{remark}

\begin{claim}
\label{a14}
If (A) then (B) where:
\mn
\begin{enumerate}
\item[$(A)$]  $(a) \quad M_1 \le_{\gk} M_2 \le_{\gk} N$
\sn
\item[${{}}$]  $(b)(\alpha) \quad M_\ell$ has cardinality $\kappa_\ell$
\sn
\item[${{}}$]  $\,\,\,\,\,\,(\beta) \quad \|N\| \ge \lambda$
\sn
\item[${{}}$]  $\,\,\,\,\,\,(\gamma) \quad \kappa_\ell \ge \kappa \ge
\LST_{\gk} + |\tau_{\gk}|$
\sn
\item[${{}}$]  $(c) \quad \Phi_1 \in \Upsilon^{\sor}_\kappa(M_1,\gk)$
\sn
\item[${{}}$]  $(d) \quad \lambda$ is strong limit and $\cf(\lambda) =
(\beth_2(\kappa_2))^+$ or just 

\hskip25pt  $(\forall \alpha < 
\cf(\lambda))(|\alpha|^{2^\kappa} < \cf(\lambda))$
\sn
\item[${{}}$]  $(e) \quad \bold x_1$ is an indirect witness for
  $(N,M_1,\lambda,\kappa,\Phi_1)$
\sn
\item[$(B)$]  there are $\Phi_2,\bold x_2$ such that:
\sn
\item[${{}}$]  $(a) \quad \Phi_2 \in \Upsilon^{\sor}_\kappa(\gk_{M_2})$
\sn
\item[${{}}$]  $(b) \quad \Phi_1 \le^1_{\kappa_2} \Phi_2$
\sn
\item[${{}}$]  $(c) \quad \bold x_2$ is an indirect witness for
  $(N,M_2,\lambda,\kappa_2,\Phi_2)$.
\end{enumerate}
\end{claim}

\begin{PROOF}{\ref{a14}}
By clause $(A)(e)$
of the assumption and the definition of indirect
witness in \ref{a12} there is $\Psi_1$ such that:
\mn
\begin{enumerate}
\item[$(*)_1$]  $(a) \quad \Psi_1 \in
\Upsilon^{\oor}_{\kappa_1}[\gk_{M_1}]$ which is standard
\sn
\item[${{}}$]  $(b) \quad \bold x_1$ is a direct witness of
$(N,M_1,\lambda,\kappa_1,\Psi_1)$
\sn
\item[${{}}$]  $(c) \quad \Phi_1 \le^4_{\kappa_1} \Psi_1$.
\end{enumerate}
\mn
By claim \ref{a9}(2) there are $\bold x_2,\Psi_2$ such that
\mn
\begin{enumerate}
\item[$(*)_2$]  $(a) \quad \Psi_2 \in
\Upsilon^{\sor}_{\kappa_2}[\gk_{M_2}]$
\sn
\item[${{}}$]  $(b) \quad \Psi_1 \le^1_{\kappa_2} \Psi_2$
\sn
\item[${{}}$]  $(c) \quad \bold x_2$ is a direct witness 
of $(N,M_2,\lambda,\kappa_2,\Psi_2)$.
\end{enumerate}
\mn
Lastly, by \ref{z30} applied to our $\Phi_1,\Psi_1,\Psi_2$ and get
$\Phi_2$ such that
\mn
\begin{enumerate}
\item[$(*)_3$]  $(a) \quad \Phi_1 \in
\Upsilon^{\sor}_\kappa[\gk_{M_2}]$ 
\sn
\item[${{}}$]  $(b) \quad \Phi_1 \le^1_\kappa \Phi_2$
\sn
\item[${{}}$]  $(c) \quad \Phi_2 \le^4_\kappa \Psi_2$.
\end{enumerate}
\mn
So we have gotten Clause (B) as promised.
\end{PROOF}

\begin{claim}
\label{a19}
If (A) + (B) then (C) where:
\mn
\begin{enumerate}
\item[$(A)$]  $(a) \quad \lambda_n \ge \LST_{\gk}$ is 
strong limit, $\cf(\lambda_n) =
(\beth_2(\LST_{\gk} + \lambda_m))^+$ if $n = m+1$
\sn
\item[${{}}$]  $(b) \quad \lambda = \sum\limits_{n} \lambda_n$ and
  $\lambda_n < \lambda_{n+1}$
\sn
\item[${{}}$]  $(c) \quad N \in K^{\gk}_\lambda$
\sn
\item[${{}}$]  $(d) \quad M_n  \le_{\gk} M_{n+1} <_{\gk} N$ and
  $\|M_n\| = \lambda_n$
\sn
\item[${{}}$]  $(e) \quad N = \cup\{M_n:n < \omega\}$
\sn
\item[$(B)$]  there is no $\Phi \in \Upsilon^{\sor}_\lambda[\gk_N]$, see
  \ref{z29}
\sn
\item[$(C)$]  for some $n$ and $\Phi$
\sn
\item[${{}}$]  $(a) \quad \Phi \in \Upsilon^{\sor}_{\lambda_n}[\gk_{M_n}]$
\sn
\item[${{}}$]  $(b) \quad$ there is an indirect witness\footnote{hence
 also a direct one; similarly in $\otimes(d)$ in the proof}
 for $(N,M_n,\lambda_{n+4},\lambda_n,\Phi_n)$ 
\sn
\item[${{}}$]  $(c) \quad$ there is no indirect witness for
$(N,M_n,\lambda_{n+5},\lambda_n,\Phi_n)$. 
\end{enumerate}
\end{claim}

\begin{remark}
1) Later we shall weaken $(A)(a)$.

\noindent
2) We may use $\Upsilon^{\sor}_\kappa[\gk_{M_n}]$ where $\lambda_0 \ge
   \kappa \ge \LST_{\gk} + |\tau_{\gk}|$ in \ref{a14} and in
   \ref{a19}, also in \ref{a22}.
\end{remark}

\begin{PROOF}{\ref{a19}}
We assume $(A) + \neg(C)$ and shall prove $\neg(B)$, this suffices. 
We try to choose $(\Phi_n,\bold x_n)$ by induction on $n$ such that:
\mn
\begin{enumerate}
\item[$\otimes$]  $(a) \quad \Phi_n \in \Upsilon^{\sor}_{\lambda_n}[\gk_{M_n}]$
\sn
\item[${{}}$]  $(b) \quad \{c_a:a \in N\} \cap \tau(\Phi_n) = \{c_a:a
\in M_n\}$
\sn
\item[${{}}$]   $(c) \quad \bold x_n$ is an indirect witness for
  $(N,M_n,\lambda_{n+4},\lambda_n,\Phi_n)$
\sn
\item[${{}}$]  $(d) \quad$ if $n=m+1$ then $\Phi_m \le^1_{\lambda_n}
 \Phi_n$.
\end{enumerate}
\mn
Now
\mn
\begin{enumerate}
\item[$(*)_1$]  if we succeed to carry the induction \then \, 
there is $\Phi \in \Upsilon^{\sor}_\lambda[\gk_N]$.
\end{enumerate}
\mn
[Why?  Note that $\Phi_n \in \Upsilon^{\sor}_{\lambda_n}[\gk_{M_n}] \subseteq
\Upsilon^{\sor}_{\lambda_n}[\gk]$ and as $\lambda_n \le \lambda$ clearly
  $\Phi_n \in \Upsilon^{\sor}_{\lambda_n}[\gk] \subseteq 
\Upsilon^{\sor}_\lambda[\gk]$ and so by \ref{z23}(2) 
there is $\Phi \in \Upsilon_\lambda[\gk]$ such that
$n < \omega \Rightarrow \Phi_n \le^1_\lambda \Phi$.  Easily $N$ is
$\le_{\gk}$-embeddable into every $\EM_{\tau(\gk)}(I,\Phi)$, in fact,
$\Phi \in \Upsilon_\lambda[\gk_N]$, contradiction to clause (B) of the
assumption.]
\mn
\begin{enumerate}
\item[$(*)_2$]   we can choose $(\bold x_n,\Phi_n)$ for $n=0$.
\end{enumerate}
\mn
[Why?  By \ref{a9}(1).]
\mn
\begin{enumerate}
\item[$(*)_3$]   if $n=m+1$ and we have chosen $(\bold x_m,\Phi_m)$
  then we can choose $(\bold x_n,\Phi_n)$.
\end{enumerate}
\mn
[Why?  If there is no indirect witness $\bold y_m$ for
  $(N,M_m,\lambda_{m+5},\lambda_m,\Phi_m)$ we have gotten 
clause (C), so \wilog \, $\bold y_m$ exists.  Now
  apply \ref{a14} with $(\bold y_n,M_m,M_n,\lambda_{n+5},\lambda_n)$ here
  standing for $(\bold x_1,M_1,M_2,\lambda,\kappa,\Phi_1)$ 
there, so we get $\bold x_n,\Phi_n$ here stand for $\bold x_2,\Phi_2$ there.]
\end{PROOF}

\begin{claim}
\label{a22}
We have $\dot I(\mu,K_{\gk}) \ge \chi$ \when \,:
\mn
\begin{enumerate}
\item[$\oplus$]  $(a) \quad \LST_{\gk} + |\tau_{\gk}| 
\le \kappa \le \chi_1 < \chi_2 < \chi_3 \le \min\{\lambda,\mu\}$
\sn
\item[${{}}$]  $(b) \quad M \le_{\gk} N$
\sn
\item[${{}}$]  $(c) \quad \|M\| \le \kappa$ and $\|N\| \ge \lambda$
\sn
\item[${{}}$]  $(d) \quad \Phi \in \Upsilon^{\sor}_\kappa[\gk_M]$
\sn
\item[${{}}$]  $(e) \quad \bold x$ is an indirect witness for
  $(N,M,\chi_2,\chi_1,\Phi)$
\sn
\item[${{}}$]  $(f) \quad$ there is no indirect witness for
$(N,M,\chi_3,\chi_1,\Phi)$
\sn
\item[${{}}$]  $(g) \quad \chi_3$ is strong limit of cofinality 
$(\beth_2(\chi_2))^+$
\sn
\item[${{}}$]  $(h) \quad \chi = |\{\theta:\theta = \beth_\theta$ and
$\theta \in [\chi_1,\chi_2]\}|$
\end{enumerate}
\end{claim}

\begin{PROOF}{\ref{a22}}
Let $\gamma_*$ be maximal such that $\beth_{\omega \cdot \gamma_*}(\chi_1)
\le \chi_2$.  Let $\Psi \in \Upsilon^{\sor}_\kappa[\gk_M]$ be 
such that $\Phi \le^4_\kappa \Psi$ and $\Psi$ has a direct witness for
$(N,M,\chi_2,\chi_1,\Psi)$ and choose such a witness $\bold x$.

Let $M_2$ be such that $M \le_{\gk} M_2 \le_{\gk} N$ and $\|M_2\| =
\beth_{\omega \cdot \gamma_*}(\chi_1) \le \chi_2$ and $\bold x$ is a direct
witness for $(M_2,M,\beth_{\omega \cdot
\gamma}(\chi_1),\chi_1,\gamma_*,\Psi)$. 

As $\chi_3$ is strong limit of cofinality $> \beth_2(\chi_2)$ there
are $\Phi_3 \in \Upsilon^{\sor}_\kappa[\gk_{M_2}]$ and $\bold y$ 
which is a direct witness for 
$(N,M_2,\chi_3,\chi_2,\Phi_3)$ and so
$\tau(\Phi_3) \backslash \{c_a:a \in M_2\}$ has cardinality $\kappa$.
For each $\gamma < \gamma_*$ there are $M_{2,\gamma},\bold x_\gamma$  
such that:
\mn
\begin{enumerate}
\item[$(*)_1$]  $(a) \quad M_{2,\gamma} \le_{\gk} M_2$
\sn
\item[${{}}$]  $(b) \quad \|M_{2,\gamma}\|$ is $\ge \beth_{\omega \cdot
\gamma}(\chi_1)$ but $< \beth_{\omega \cdot \gamma + \omega}(\chi_1)$;
  can get even 

\hskip25pt $\|M_{2,\gamma}\| = \beth_{\omega \cdot \gamma}(\chi_1)$
\sn
\item[${{}}$]  $(c) \quad \bold x_\gamma$ is a direct witness for
$(M_{2,\gamma},M,\beth_{\omega \cdot
  \gamma}(\chi_1),\chi_1,\gamma,\Psi)$.
\end{enumerate}
\mn
[Why?  Try by induction on $k$ to choose $\eta_k \in \cT_{\bold x}$
such that $\ell g(\eta_k) = 2k +1,\eta_k(2k) \ge \beth_{\omega \cdot
\gamma}(\chi_1)$ and $\ell < k \Rightarrow \eta_\ell \triangleleft
\eta_k$.  For $k=0$, clearly $\eta_k = \langle \rangle$ is O.K., and
as $\eta_\ell(2 \ell) > \eta_{\ell +1}(2 \ell +2)$, necessarily for
some $k$ we have $\eta_k$ but cannot choose $\eta_{k+1}$; let $A_\gamma
= \cup\{\Rang(h_{\eta,u}):\eta_k \trianglelefteq \eta \in \cT_{\bold
y}$ and $u \in [I_\eta]^{\lfloor \ell g(\eta)/2 \rfloor}\}$ so
$A_\gamma \subseteq M$ has cardinality $\eta_k(2k) \in [\beth_{\omega
\cdot \gamma}(\chi_1),\beth_{\omega \cdot \gamma + \omega}(\chi_1)$.
\Wilog \, if $N_* = \EM(\emptyset,\Phi_3)$ is a standard (i.e. $M = N_*
\rest \tau_{M_2})$ then $A_\gamma$ is closed under the
functions of $N_* \rest \tau'_{\Phi_3}$.  Let $M_{2,\gamma} = M_2
\rest A_\gamma$; it is $\le_{\gk} M$ and it 
satisfies clauses (a),(b) and include $A_*$.  Then we can easily find
$\bold x_\gamma$ as required in clause (c).]

Next we can find $\bold y_\gamma,\Phi_{3,\gamma}$ such that
\mn
\begin{enumerate}
\item[$(*)_2$]  $(a) \quad \bold y_\gamma$ is a direct witness of
$(N,M_{2,\gamma},\chi_3,\|M_{2,\gamma}\|,\Phi_{3,\gamma})$
\sn
\item[${{}}$]  $(b) \quad \Phi_{3,\gamma} \in
\Upsilon^{\sor}_\kappa[M_{2,\gamma},\gk]$.
 \end{enumerate}
\mn
[Why?  Recall $\tau(\Phi_3) \backslash \{c_a:a \in M_2\}$ has
  cardinality $\kappa$.  Let $\tau_{2,\gamma} = \tau(\Phi_3)
  \backslash \{c_a:a \in M_2 \backslash M_{2,\gamma}\}$ so has
  cardinality $\|M_{2,\gamma}\|$, let $\Phi_{3,\gamma} = \Phi_3 \rest
  \tau_{2,\gamma}$, is as required in $(*)_2(k)$.  As for $\bold
  y_\gamma$ we derived it form $\bold y$.]

Now let $I=I_\mu$ be a linear order of cardinality $\mu$ as 
required in \ref{z9}.

Lastly, let $N_\gamma = \EM_{\tau(\gk)}(\mu,\Phi_{3,\gamma})$ be
standard hence $M_{2,\gamma} \le_{\gk} N_\gamma \in K^{\gk}_\mu$.

We choose $\partial_i$ by induction on $i$ such that: if $i=0$ then
$\partial_i = \chi_1$, if $i$ is limit then $\partial_i = \cup\{\partial_j:j
< i\}$ and if $i=j+1$ then $\partial_i = \beth_{\beth_2(\partial_j)^+}$
when it is $\le \chi_2$ and undefined otherwise.  Let $\partial_i$ 
be defined iff $i<i(*)$ and let $\Theta = \{\partial_{i+1}:i +1 < i(*)\}$.  Now
$|\Theta| \ge \chi$ so it suffices to prove that
$\langle N_\theta:\theta \in \Theta\rangle$ are pairwise
non-isomorphic.

So toward contradiction assume
\mn
\begin{enumerate}
\item[$(*)_3$]  $\theta_1 < \theta_2$ are from $\Theta$ and $\pi$ is an
isomorphism from $N_{\theta_2}$ onto $N_{\theta_1}$.
\end{enumerate}
\mn
We can find $M_* \le_{\gk} N_{\theta_1}$ such that $\|M_*\| =
\theta_2$ and $M \cup M_{2,\theta_1} \cup \pi(M_{2,\theta_2}) \subseteq M_*$
and \wilog \, we can find 
$I_* \subseteq \mu$ of cardinality $\theta_2$ such that
$M_* = \EM_{\tau(\gk)}(I_*,\Phi_{3,\theta_1})$.

Let $I^*_1 \subseteq I_*$ be of cardinality $\theta_1$ such that
$M_{2,\theta_1} \cup \pi(M) \subseteq N'_{\theta_1} :=
\EM_{\tau(\gk)}(I^*_1,\Phi_{3,\theta_1})$ and let $N'_{\theta_2} =
\pi^{-1}(N'_{\theta_1})$.  By \ref{a9}(2) we can find $\Psi' \in
\Upsilon^{\sor}_\kappa(N'_{\theta_2},\gk)$ and $\bold x_{\theta_2}$ a witness
  for $(M_{2,\theta_2},N'_{\theta_2},\theta_2,\kappa,\Psi')$ such that
$\Psi \le^4_\kappa \Psi'$ and $\bold x_{\theta_2} \le \bold
  x'_{\theta_2}$ where $\theta_2 = \beth_{\omega \cdot \gamma_2}(\chi_1)$.

Now clearly $N'_{\theta_1},\Psi,\pi(\Psi'),\pi(\bold x'_{\theta_2})$
satisfies the parallel statements in $N_{\theta_1}$.  By
\ref{a13}(B)(a) and the choice of $I_\mu$ there is a witness for
$(N_{\theta_1},N'_{\theta_1},\chi_3,\kappa,\pi(\Psi'))$, hence applying
  $\pi^{-1}$ there is a witness $\bold x''_{\theta_2}$ for
  $(N_{\theta_1},N'_{\theta_1},\chi_3,\kappa,\Psi')$.

Hence by \ref{a13}(B)(b), $\Psi' \le^3_\kappa \Phi_{3,\theta_2}$ but
together $\Phi \le^4_\kappa \Psi \le^4_\kappa \Psi' \le^3_\kappa
\Phi_{3,\theta_2}$ hence $\Phi \le^3_{\theta_2} \Phi_{3,\theta_2}$ by
\ref{z26}(1) so by \ref{z26}(2), the last clause, there is
$\Phi'_{3,\theta_2} \in \Phi_{3,\theta_2}/\bold E^{\ai}_{\theta_2}$
such that $\Phi \le^4_{\theta_2} \Phi'_{3,\theta_2}$.  
But as $\Phi_{3,\theta_2}$
has a $(N,M_{2,\theta_2},\chi_3,\theta_2)$ witness by \ref{a7}(3) also
$\Phi'_{3,\theta_2}$ has hence $\Phi$ has an indirect witness for
$(N,M,\chi_3,\kappa)$, contradiction.
\end{PROOF} 

\begin{conclusion}
\label{a25}
Assume $\cf(\lambda) = \aleph_0$ and $\lambda = \beth_{1,\lambda}$.

\noindent
1) If $\lambda > \dot I(K_{\gk})$ then 
$M \in K^{\gk}_\lambda \Rightarrow \Upsilon^{\sor}_\lambda[\gk_M] 
\ne \emptyset$.  

\noindent
2) If $\mu \ge \lambda > \dot I(\mu,K_{\gk})$ \then \, $M \in
K^{\gk}_\lambda \Rightarrow \Upsilon^{\sor}_\lambda[\gk_M] \ne \emptyset$.

Moreover, at least one of the following holds:
\mn
\begin{enumerate}
\item[$(a)$]  for some $\chi_1 < \lambda$ if $\chi_1 < \chi_2 =
\beth_{2,\delta} \le \min\{\lambda,\mu\}$ 
then $|\delta| \le \dot I(\mu,K_{\gk})$
\sn
\item[$(b)$]  $\Upsilon^{\sor}_\lambda[\gk_M] \ne \emptyset$ for
every $M \in K^{\gk}_\lambda$.
\end{enumerate} 
\end{conclusion}

\begin{theorem}
\label{a28}
The result from the abstract holds, that is, for every a.e.c. $\gk$
for some closed unbounded class $\bold C$ of cardinals we have 
(a) or (b) where
\mn
\begin{enumerate}
\item[$(a)$]  for every $\lambda \in \bold C$ of cofinality $\aleph_0,
\dot I(\lambda,K) \ge \lambda$ 
\sn
\item[$(b)$]  for every $\lambda \in \bold C$ of cofinality $\aleph_0$ and 
$M \in K_\lambda$, for every cardinal $\kappa \ge \lambda$ there is
$N_\kappa$ of cardinality $\kappa$ extending $M$ (in the sense of our
a.e.c.).
\end{enumerate}
\end{theorem}

\begin{PROOF}{\ref{a28}}
Let $\Theta = \{\mu:\mu = \beth_{2,\delta}$ and $|\delta| > \dot
I(\mu,K_{\gk})$ for some limit ordinal $\delta\}$.
\medskip

\noindent
\underline{Case 1}:  $\Theta$ is an unbounded class of cardinals.

So $\bold C = \{\mu:\mu = \sup(\mu \cap \Theta)\}$ is a closed
unbounded class of cardinals.  Easily $\mu \in \bold C \Rightarrow \mu
= \beth_{1,\mu}$ and by \ref{a25} + \ref{z29} for every $\mu \in \bold
C$, clause (b) of \ref{a28} holds.
\medskip

\noindent
\underline{Case 2}:  $\Theta$ is a bounded class of cardinals.

So by the definition of $\Theta,\bold C = \{\mu:\mu > \sup(\Theta),\mu
= \beth_{2,\mu}\}$ is as required.
\end{PROOF}

\noindent
Also
\begin{theorem}
\label{a30}
For every $\aec$ \, $\gk$ one of the following holds:
\mn
\begin{enumerate}
\item[$(a)$]  for some $\chi$ we have $\chi < \mu = \beth_{2,\mu}
  \Rightarrow \dot I(\mu,K_{\gk}) \ge \mu$ and $\chi < \mu =
  \beth_{1,\omega \cdot \gamma} \Rightarrow \dot I(\mu,K_{\gk}) \ge
  |\gamma|$
\sn
\item[$(b)$]  for some closed unbounded class $\cC$ of cardinals we
  have $\cf(\lambda) = \aleph_0 \wedge \lambda \in \bold C \wedge M
 \in K^{\gk}_\lambda \Rightarrow \Upsilon^{\sor}[M,\gk] \ne
  \emptyset$.
\end{enumerate}
\end{theorem}

\begin{PROOF}{\ref{a30}}
Similarly to \ref{a28}, using Fodor lemma for classes of cardinals.
\end{PROOF}
\newpage

\section {Concluding Remarks} \label{Concluding}

\begin{definition}
\label{n2}
1) For an ordinal $\gamma,\tau$-models $M_1,M_2$ and cardinal
   $\lambda$ we define a game $\Game =
   \Game_{\theta,\gamma}(M_1,M_2)$.  A play lasts less than $\omega$
   models is defined as in \cite[2.1]{Sh:797}.
\end{definition}

\begin{claim}
\label{n4}
1) Assume $\cf(\lambda) = \aleph_0$ and $M_1,M_2$ are $\tau$-models of
   cardinality $\lambda$.  If the isomorphic player wins in
   $\Game_{\lambda,\gamma}(M_1,M_2)$ for every $\gamma$ or just
   $\gamma < (2^{< \lambda})^+$ \then \, $M_1,M_2$ are isomorphic.

\noindent
1A) If above $\lambda$ is strong limit  \then \, ``$(2^{< \lambda})^+ =
\lambda^+$".

\noindent
2) Assume $\lambda$ is strong limit of cofinality $K =
 K_{\gk}$ and $|\tau_{\gk}| + \LST_{\gk} \le \lambda$ and $K = \{M \rest
   \tau:M \models \psi\}$ for some $\psi \in
   \bbL_{\lambda^+,\aleph_0}$.

If $\dot I(\lambda,K) \le \lambda$ \then \, for every $M_1 \in K$
there is $M_2 \in K_{\le \lambda}$ such that the isomorphic player
wins in $\Game_{\lambda,\gamma}(M_1,M_2)$ for every $\lambda$.
\end{claim}

\begin{conjecture}
\label{n8}
For every a.e.c. $\gk$ letting $\kappa = \LST_{\gk} + |\tau_{\gk}|$,
at least one of the following occurs:
\mn
\begin{enumerate}
\item[$(a)$]  if $\lambda = \beth_{1,\lambda} > \kappa$ and $\cf(\lambda) =
  \aleph_0$, \then \, $\Upsilon^{\sor}_\kappa[M,\gk] \ne \emptyset$
\sn
\item[$(b)$]  if $\lambda = \beth_{1,\lambda} > \kappa$ and $\cf(\lambda) =
\aleph_0$, \then \, $\dot I(\lambda,K_{\gk}) = 2^\lambda$.
\end{enumerate}
\end{conjecture}
\newpage

\bibliographystyle{alphacolon}
\bibliography{lista,listb,listx,listf,liste,listz}

\end{document}